\documentclass{article}   % or amsart
\usepackage{theorem}    % for the \theoremstyle, cf. Samarin
\usepackage{version}

\usepackage{amsmath}
\usepackage{amsfonts}
\usepackage{color}      % for the xfig figures
\usepackage{graphicx}   % for the ps pictures
\usepackage[all]{xy}
\numberwithin{equation}{subsection}
\theorembodyfont{\rmfamily}

\newtheorem{thing}{}[subsection]
\theoremstyle{definition}  % cf. Samarin, Gratzer

\newcommand{\dual}{^{\vee}}

\newcommand{\jk}{{J_{\xi}}}
\newcommand{\ok}{{\OO_{\xi}}}
\newcommand{\ox}{{\OO_{X}}}
\newcommand{\oc}{{\OO_{C}}}

\newcommand{\os}{{\OO_{S}}}

\newcommand{\OO}{{\cal O}} % Structure sheaf
\newcommand{\PP}{{\mathbb P}} % projective space
\renewcommand{\AA}{{\mathbb A}} % affine space
\newcommand{\Spec}{\mathop{\rm Spec}}

\newcommand{\pic}{{\rm Pic}}
\newcommand{\Num}{\mathop{\rm Num}}
\newcommand{\hilb}{{\rm Hilb}}
\newcommand{\Hilb}{{\rm Hilb}}
\newcommand{\Pic}{\mathop{\rm Pic}\nolimits}
\newcommand{\vdim}{\mathop{\rm v.dim}\nolimits}
\newcommand{\triv}[1]{#1 \tensor \OO_X}   % V tensor \OO

\newcommand{\ZZ}{\mathbb Z}

\newcommand{\RR}{\mathbb R}
\newcommand{\rk}{\mathop{\rm rk}}
\newcommand{\rank}{\mathop{\rm rank}}
\newcommand{\codim}{\mathop{\rm codim}\nolimits}

\newcommand{\Ext}{\mathop{\rm Ext}\nolimits}  
\newcommand{\Hom}{\mathop{\rm Hom}\nolimits}

\newcommand{\sheafHom}{\mathcal H \mathit{om}}
\newcommand{\coker}{\mathop{\rm coker}}

\newcommand{\tensor}{\mathop{\otimes}}

\newcommand{\ch}{\mathop{\rm ch}}  
\newcommand{\Gr}{\mathop{\rm Gr}}

\newcommand{\intersect}{\cap}

\newcommand{\isom}{\stackrel{\sim}{\to}}

\newcommand{\strelka}[1]{\stackrel{#1}{\to}} % arrow s nadpis'yu

\newcommand{\triple}[3]{
                         0 \to {#1} \to {#2} \to {#3} \to 0 
                       }

\newcommand{\xyArrow}[2]{\ar@{}[#2]|-*{#1}}
\newcommand{\xySimeq}{\ar@{-}^*@{~}}      % \simeq - strelka 
\newcommand{\xyRavno}{\ar@{=}}
\newcommand{\xyMaps}{\ar@{|->}}           % element perehodit v element
\newcommand{\xyKrivayastrelka}{\ar@{~>}}  % krivaya strelka
\newcommand{\xyMono}{ \ar@{^{(}->}}       %vlozhenie
\newcommand{\xyMonom}{\ar@{_{(}->}}       %tozhe vlozhenie
\newcommand{\xyEpi}{\ar@{->>}}            %epimorfizm
\newcommand{\xyCherta}{\ar@{-}}           %pryamaya liniya
\newcommand{\xySubsetSmall}{\ar@{}|(0.60){\subset}} % malen'kaya strelka subset
\newcommand{\xySubset}{\ar@{}|(0.60)*{\subset}} % subset
\newcommand{\xySubsetA}{\ar@{}|*{\subset}[r]} % subset, raspolozhennyj
\newcommand{\xyCap}{\ar@{}|(0.60)*{\cap}} % cap (vlozhenie vniz)
\newcommand{\xyCapA}{\ar@{}|*{\cap}[d]}   %      vlozhenie vniz, raspolozheno
\newcommand{\xyCup}{\ar@{}|*{\cup}[u]}    %      vlozhenie vverh 

\newcommand{\smallXyMatrix}{\xymatrix@R=10pt@C=10pt}

\newcommand{\correspondence}[5]
{ 
   \smallXyMatrix
   {
        {} & {{#1}}\ar[ld]_-{{#4}}\ar[rd]^-{{#5}} & {} \\
        {{#2}} & {} & {{#3}}
   }
}   
    
\newcommand{\ep}{\epsilon}
\newcommand{\cA}{{\cal A}}
\newcommand{\smirr}{^{s}}     % sm. irr. (smooth irreducible)
\newcommand{\br}{\bar{r}}
\begin{document}

\title{On the Brill-Noether theory for K3 surfaces}
\author{Maxim Leyenson}
\date{}
\maketitle
%-------------------------------------------------
\begin{abstract}
    Let $(S,H)$ be a polarized K3 surface. We define Brill-Noether
filtration on moduli spaces of vector bundles on $S$. Assume that
$(c_1(E),H) > 0$ for a sheaf $E$ in the moduli space. We give a
formula for the expected dimension of the Brill-Noether
subschemes. Following the classical theory for curves, we give a
notion of Brill-Noether generic K3 surfaces.

   Studying correspondences between moduli spaces of coherent sheaves
of different ranks on $S$, we prove our main theorem: polarized K3
surface which is generic in sense of moduli is also generic in sense
of Brill-Noether theory (here $H$ is the positive generator of the
Picard group of $S$). In case of algebraic curves such a theorem,
proved by Griffiths and Harris and, independently, by Lazarsfeld, is
sometimes called ``the strong theorem of the Brill-Noether theory''.

We finish by considering a number of projective examples. In
particular, we construct explicitly Brill-Noether special K3 surfaces
of genus 5 and 6 and show the relation with the theory of
Brill-Noether special curves.
\end{abstract}
%-------------------------------------------------

\tableofcontents
\subsection{Introduction}
\subsubsection{Classical Brill-Noether theory.}
    We start by recalling some results of classical Brill-Noether
theory. Let $C$ be a smooth projective curve of genus $g$ over an
algebraically closed field $k$ of characteristics 0; consider the
Picard variety $\Pic^d(C)$ parameterizing degree $d$ line bundles on
$C$. For a line bundle $L$ on $C$ one can define cohomology groups
with coefficients in $L$, $H^0(C,L)$ and $H^1(C,L)$; these are vector
spaces over the base field $k$. The vector space $H^0(C,L)$ can be
interpreted as the space of global sections of $L$. Let us denote by
$h^i(C,L)$ the dimension of the vector space $H^i(C,L)$ over $k$, $i =
0, 1$. These numbers are related by the Riemann-Roch formula
$$
     h^0(C,L)-h^1(C,L) = \chi(L) = d+1-g,
$$
thus, given $h^0(C,L)$, we can compute $h^1(C,L)$ from the
formula above.

   Therefore, the Picard variety $\Pic^d(C)$ has a natural filtration
given by the value of $h^0(C,L)$, 
$$
   \Pic^d(C) \supset W_d^0(C) \supset W_d^1(C) \supset \dots,
$$
where $W_d^r(C)$ is the set of isomorphism classes of the degree $d$
line bundles $L$ which satisfy the condition $h^0(L) \ge r+1$. (The shift 
by +1 appears for historical reasons: if $h^0(C,L) = r+1$, then $r$ is the 
dimension of the linear system $|L|$).

  One can prove that each $W_d^r(C)$ is a set of points of a closed
subvariety in the Picard variety $\Pic^d(C)$, and it is classically
known that many questions of the projective geometry of $C$ can be
reformulated in terms of the geometry of these subvarieties, named
``the Brill-Noether loci'' after Brill and Noether who did the first
essential study of them in \cite{brill-noether} in 1873. (The notation
$W_d^r$ itself goes back to Brill and Noether.) For example, the
variety $W_2^1(C)$ is not empty iff the curve $C$ is hyperelliptic,
and if $C$ can be represented as a plane curve of degree $d$, then
$W_d^2(C)$ is not empty.

    What is known about the geometry of the Brill-Noether loci for a
curve $C$?  We formulate only a few results. First, one can construct
two vector bundles, $E$ and $F$, on a Picard variety $\Pic^d(C)$, and
a map $\sigma: E \to F$ in such a way that $W_d^r$ is equal to the
degeneration subvariety of $D_k(\sigma)$ for a certain $k$, where
$$
   D_k(\sigma) = \{ x \in \Pic^d(C): \rank \sigma(x) \le k \}
$$
(see, e.g., \cite{fulton}, Chapter 14.)
%given by
%the condition that the rank of $\sigma$ drops by .... .  
In such a situation one can define the ``expected'' (or ``virtual'')
dimension of the variety $W_d^r$ in such a way that if the map
$\sigma$ is ``generic'' in an appropriate sense, then the dimension of
the degeneration variety $W_d^r$ is equal to its expected
dimension. The expected dimension of $W_d^r(C)$ depends only on the
numbers $g$, $d$ and $r$, but not on the curve $C$ itself, and is
classically denoted as $\rho(r,g,d)$. One has
$$   
   \rho(r,g,d) = g - (r+1)(g-d+r)
$$
and $\rho(r,g,d)$ is called ``the Brill-Noether number''.

If, for a given curve $C$ and numbers $r$ and $d$, the variety
$W_d^r(C)$ has the ``expected'' dimension $\rho(r,g,d)$, then much is
known about the geometry of $W_d^r(C)$, for example, one can compute
its cohomology class (in the case $k = \mathbb{C}$) in the cohomology
ring of the Picard variety (A. Poincare), its normal invariants, and
so on. (See \cite{ACGH} for many more details).

  The curve $C$ is said to be ``Brill-Noether general'' if 
\begin{itemize}
   \item[(a)] whenever $\rho(r,g,d) < 0$ (the expected dimension of
$W_d^r(C)$ is negative), the variety $W_d^r(C)$ is empty, and
   \item[(b)] whenever $\rho(r,g,d) \ge 0$, the variety $W_d^r(C)$ is
not empty and $\dim W_d^r(C) = \rho(r,g,d)$ (i.e., $W_d^r(C)$ is of
expected dimension).
\end{itemize}

    It was conjectured by Brill and Noether that a generic curve $C$
in the moduli space of curves of genus $g$ is (what we call now)
Brill-Noether general.  The first proof was given by Griffiths and
Harris in 1980 (\cite{griffiths-harris-1980}) and used a degeneration
argument. %
% **
%
It was also observed (Miles Reid, Tyurin, Donagi-Morrison, Lazarsfeld)
% add the references !!
that if the curve $C$ can be realized as a hyperplane section of some
K3 surface $S$, then the geometry of the Brill-Noether loci $W_d^r(C)$
is closely related to the geometry of the surface $S$. Lazarsfeld
proved that if $C$ is a curve which can be embedded into a K3 surface
$S$ in such a way that the linear system $|C|$ does not contain
non-reduced or reducible curves, then whenever $\rho(r,g,d) < 0$ the
variety $W_d^r(C)$ is empty.

\subsubsection{Generalizations.}

 There are two most natural ways to generalize the classical
Brill-Noether theory: first, instead of studying line bundles one can
consider coherent sheaves of arbitrary rank $r$. (Unfortunately, the
dimension of the linear system $L$, as $r$ in the notation $W_d^r$,
and the rank of a vector bundle are classically denoted by the same
letter $r$.) Second, instead of studying curves, one can consider
varieties $X$ of dimension greater than one. One can go in both
directions simultaneously, fixing a variety $X$ equipped with an ample
divisor class $H$ and considering the moduli space $M = M_H(r,c)$ of
$H$ - stable vector bundles (or coherent sheaves) of rank $r$ and
Chern class $c$, and studying the subschemes in $M$ defined by the
conditions $\{ h^i(X,E) \ge k_i \}$.

%------------ **  rewrite' better
 
    Going in the first direction, Newstead and others studied the
Brill-Noether loci in the moduli spaces of stable rank $r$ vector
bundles on curves (their results are published in a sequence of papers
in 1994-2004). Mumford's study of special 0-cycles, done in the
Chapters 19-21 of ~\cite{mumford} in order to construct the Picard
scheme of an algebraic surface $S$, can be interpreted as a proof of
non-emptiness of certain Brill-Noether varieties in the moduli space
of rank 1 torsion free sheaves on $S$. Goettsche and Hirschowitz
studied the Brill-Noether loci in the moduli spaces of stable vector
bundles on the projective plane $\mathbb{P}^2$. %
%
% Goettsche, Lothar; Hirschowitz, Andre
%   Weak Brill-Noether for vector bundles on the projective plane.    
%   Lecture Notes in Pure and Appl. Math., 200, Dekker, New York, 1998
%
 It is known that if $S$ is an algebraic surface over $k =
\mathbb{C}$, then the geometry of the Brill-Noether loci in the moduli
spaces of vector bundles on $S$ is connected with the invariants of
smooth structure of $S$.%
% considered as a smooth 4-dimensional manifold. 
(For example, Ellingsrud, Stromme and Le Potier studied the
Brill-Noether stratification in the moduli space of stable vector
bundles on $\PP^2$ in order to compute some Donaldson invariants of
$\PP^2$).

%---------------------------------------------------------------

%
\subsubsection{The case of K3 surfaces.}

    We consider the Brill-Noether filtration in the moduli space of
stable vector bundles on K3 surfaces.

    In general, a vector bundle $E$ on an algebraic surface $S$ has
three Betti numbers, $h^0(S,E)$, $h^1(S,E)$, and $h^2(S,E)$, which are
related by the Riemann-Roch theorem
$$
    h^0(S,E) - h^1(S,E) + h^2(S,E) = \chi(S,E),
$$
and therefore one naturally has to consider a bifiltration in the
moduli spaces of stable vector bundles on $S$ given by the condition
$h^0(S,E) \ge k, h^1(S,E) \ge l$.  However, if $S$ is a K3 surface,
one can prove that the second cohomology group $H^2(S,E)$ vanish for
all stable vector bundles $E$ such that $(c_1(E),H) > 0$. It follows
that one has to consider only two Betti numbers which are subject to
one relation, and therefore the situation in the case of K3 is simpler
and it is enough to study the filtration defined by the dimension of
the spaces of global sections.

    Let $S$ be a K3 surface. Let $a \in \Pic(S)$, $d \in \ZZ$, and let
$M = M_H(r,a,d)$ be the moduli space of $H$-stable vector bundles $E$
on $S$ which satisfy $\rk(E) = r$, $c_1(E) = a$ and $c_2(E) = d$. We
define the Brill-Noether filtration on $M$ as 
$$
    BN_k(M)        \stackrel{def}{=}
    BN_k(r,a,d)(S) \stackrel{def}{=}
    \{[E] \in M_H(r,a,d): h^0(S,E) \ge k \}.
$$

  We realize $BN_k(M)$ as the degeneration variety of a morphism of
vector bundles $\sigma: E \to F$ on the moduli space $M$. We then
compute the expected dimension of the Brill-Noether loci $BN_k(M)$,
which we denote as $\vdim BN_k(M)$ (``v'' stands for ``virtual'').

   It is natural to ask, then, whether, at least for a ``general'' K3
surface $S$, the numeric condition $\vdim BN_k(M) < 0$ implies that
the variety $BN_k(M)$ is empty, and whether $\vdim BN_k(M) \ge 0$
implies that $BN_k(M)$ is not empty and is of the expected dimension.

\begin{comment}
 We say that a triple $(S,H,a)$, where $S$ is a K3-surface, $H \in
\Pic S$, $a \in \Pic(S)$, and $H$ is an ample bundle, is Brill-Noether
general if

(a) for any given r,d and k, such that
$v.dim BN_k(M) < 0$, $M = M_H(r,a,d)$, then $BN_k(M)$ is empty, 

(b) if, for a given r,d, and k, $v.dim BN_k(M) > 0$, then $BN_k(M)$ is
not empty and is of the expected dimension.  %

    As in the classical Brill-Noether theory, it is natural to ask
whether a generic algebraic K3 surface is generic in sense of the
Brill-Noether theory, i.e.

. (For a generic algebraic K3 surface S, one has $\Pic S \simeq Z$,
and therefore we can choose $H = h$ and $a = a_0 h$, where $h$ is a
positive generator of $\Pic S$, and $a_0 \in \ZZ$.

   It turns out that in many cases K3 surfaces have a plenty of
Brill-Noether loci which expected dimension is greater than 0 but
which are empty, but we modify the question by defining a region 
\end{comment}

%   should be rewritten

    We give a partial answer to this question. To summarize our
results shortly, 
\begin{itemize}
   \item[(1)] we give examples of situations when the expected
dimension of the Brill-Noether subscheme $BN_k(v)$ is greater than
zero, but the corresponding subscheme is empty, in contrast with the
classical case, and define a ``good'' geographic region in the space
of numerical parameters $(g,d,r,k)$ where we expect the Brill-Noether
filtration to behave in the expected way for a generic K3 surface;
\item[(2)] we prove that if $S$ is a K3 surface with $\Pic S \simeq
\ZZ a$, where $a$ is the positive generator of $\Pic S$, and the
generic curve in the linear system $|a|$ is Brill-Noether general in
the classical sense, then the Brill-Noether filtration on the moduli
space $M_H(r,a,d)$ has the expected behavior in the good geographic
region, and consists of empty varieties outside of the good region, i.e.,
in this case we have a full answer to the question we study.
 A theorem of Lazarsfeld implies that the second condition above is
extra, i.e., that for a K3-surface $S$ with $\Pic S \simeq \ZZ a$,
where $a$ is the positive generator of $\Pic S$, the generic curve $C$
in the linear system $|a|$ is Brill-Noether general. Thus a generic
polarized K3 surface is Brill-Noether general in the ``good'' region.
\end{itemize}

  We also give some examples of Brill-Noether-special K3 surfaces
(with rank of $\Pic S > 1$), and give some applications to the
birational geometry of the moduli spaces of vector bundles on K3
surfaces.

\subsubsection{Techniques.}
    Below I summarize the techniques we use.

    First, we extend the problem to include certain coherent sheaves
of rank 0 and 1, thus extending our picture to the moduli spaces of
all the ranks $r \ge 0$. For rank 0, studying the corresponding
Brill-Noether filtration on the moduli space $M(0,a,d)$ is equivalent
to studying the classical Brill-Noether theory for curves on $S$ in
the linear system $|a|$, and for sheaves of rank 1 studying of the
Brill-Noether filtration on $M(1,a,d)$ is equivalent to the study of
the 0-cycles on $S$ which fail to impose independent conditions on
curves in the linear system $|a|$.

    Then, continuing the line of thought of Tyurin in \cite{tyurin},
we define some system of correspondences which relate Brill-Noether
varieties in the moduli spaces of sheaves of various ranks. We can use
this system of correspondences to obtain information about the
geometry of the Brill-Noether variety $BN_k(r,a,d)$ for a given $k$
and $r$, provided that we can guarantee that for a certain $r'$,
depending on $r,k,a$ and $d$,
\begin{itemize}
   \item[(1)] for a generic vector bundle $[E] \in M_H(r',a,d)$ one
has $h^1(S,E) = 0$, and
   \item[(2)] generic vector bundle in $M_H(r',a,d)$ is globally
generated.
\end{itemize}

     However, this seem to be a hard question on its own. In our study
of the Brill-Noether theory on $S$ we apply the correspondence twice,
the first time to translate some information about the special linear
systems on curves on the surface $S$ into information about the
behavior of the generic vector bundles $[E] \in M_H(r',a,d)$ for
various ranks $r'$, thus reversing the line of though of Lazarsfeld
and answering, under certain conditions, the question about the
behavior of the generic vector bundles posed above, and the second
time to translate the information about the generic vector bundles in
the moduli spaces $M_H(r',a,d)$ into the information about the given
Brill-Noether variety $BN_k(r,a,d)$.

%
\begin{comment}
    In order to apply the first step, we need an existence of a big
number of globally generated line bundles on curves on the surface
S. In particular, we prove that if the linear system .... contains at
least one Brill-Noether general curve, then the Martens-Mumford
theorem implies that this curve contains enough globally generated
line bundles for the Step one to work.

     It is generically considered that in the series of papers
.... Lazarsfeld proved that if $\Pic S = Z$, (which is the case of the
most general algebraic K3 surface), then a generic curve in the linear
system [which is a generator of the Picard group of S] is
Brill-Noether general, therefore allowing us to complete Step
1. However, I was not able neither to understand some details of the
Lazarsfeld's proof, nor to proof the result myself, nor to find a
counterexample. Therefore, I first prove a slightly weaker statement
and then re-state it admitted that the Lazarsfeld's results are true.
\end{comment}
%
\subsubsection{Structure of the paper.}
The proof of some technical results will appear in a separate
publication \cite{bn-K3-part2}. The proof of the main theorem is given
completely in this paper.
\subsubsection{Acknowledgments.}
My interest in the Brill-Noether theory was prompted by a conversation
with prof. M.S. Narasimhan at the International Center for Theoretical
physics (ICTP), Trieste, Italy in 1999. I would like to thank
prof. M.S. Narasimhan and ICTP for hospitality.

The proof of the main theorem is based on the ideas of Robert
Lazarsfeld and Andrey Tyurin. This paper is a result of my attempt to
understand the ``Cycles, curves and vector bundles on algebraic
surfaces'' by Andrey Tyurin (~\cite{tyurin}).

I am thankful do Dima Arinkin for reading the preliminary version of
the manuscript. 

I am especially thankful to Andrey Levin of the Institute for
Theoretical and Experimental Physics (ITEP, Moscow) for numerous
discussions, help, encouragement and psychological support.

The proof of the technical results mentioned in
\label{technical} will hopefully appear in the second part of the paper.

\subsection{Preliminaries: moduli spaces
             of acceptable sheaves on a K3 surface.}
Let $S$ be a K3 surface over an algebraically closed field $k$ of
characteristic 0. Let $K(S) = K^0(S)$, $K'(S) = \ZZ \oplus \Num S
\oplus \ZZ$, and $K'' = \ZZ^3$. Consider the maps
$$
   K(S) \stackrel{\ch}{\isom} K'(S) \stackrel{\epsilon}{\to} K'',
$$
% K'' does not depend on S
%
where the Chern character map $\ch$ takes a class of a locally free
sheaf $E$ to its Chern character $\ch(E) = (\rk E, \ch_1(E),
\ch_2(E))$, and $\ep(r,v_1,v_2)=(r,g,d)$, where $g=\frac{v_1^2+2}{2}$,
and $d = \frac{v_1^2-2v_2}{2}$. Note that $\ch_2(E) = (c_1(E)^2 -
2c_2(E))/2 \in \ZZ$, since $\Num S$ is an even lattice, so the map
$\ch$ is correctly defined. Note also that for a locally free sheaf
$E$ one has $\ep(v([E])) = (\rk E,g,d)$, where $g$ is the arithmetic
genus of a curve in the linear system $|c_1(E)|$ in the case it is not
empty, and $d=c_2(E)$. Note that linear equivalence on $S$ coincides
with numerical equivalence, since $H^1(S,\OO_S) = 0$.

   If $a \in \Num S$, we denote by $|a|$ the linear system $\PP
H^0(S,L_a)\dual$, where $L_a$ is the unique invertible sheaf on $S$
with $c_1(L_a) = a$, and by $|a|\smirr$ the open subscheme in $|a|$
parametrizing smooth irreducible curves.

We say that a coherent sheaf $F$ on $S$ is {\it acceptable}
\footnote{This is not a standard terminology} if it satisfies one of
the following conditions:\\
\begin{enumerate}
  \item $\rk F \ge 2$ and $F$ is locally free; 
  \item $\rk F = 1$ and $F$ is torsion free; 
  \item $\rk F = 0$ and $F$ is isomorphic to the direct image of an invertible
sheaf on a smooth irreducible curve on $S$.  
\end{enumerate}
We give some justification to this definition later. 

Let $H \in Num(S)$ be a class of an ample divisor, let $v= (r,v_1,v_2)
\in K'(S)$ such that $(v_1,H) > 0$, and let $M(v)=M_H(v)$ be the
moduli space of $H$-stable acceptable coherent sheaves on $S$ with
Chern character $v$. (This moduli space can be constructed as an open
subscheme in the moduli space of stable coherent sheaves. For the
construction of the latter see, for example, \cite{maruyama}.)

 Let $\ep(v)=(r,g,d)$. We define $$\beta = v_2 = g-1-d$$ and $$\alpha
= - \beta = d+1-g$$

\begin{thing} {\bf Definition.} 
$$
    \rho(r,g,d):= g - (r+1)(r+g-d)
$$
\end{thing}

Note that if $d < 0$, then $\rho(r,g,d) < 0$.

For any $r \in \ZZ$ we let $\br = r-1$.

\begin{thing} {\bf Theorem (Shigeru Mukai, \cite{mukai-84,mukai-84a})}. If
$r \ge 2$ and $\rho(\br,g,d) \ge 0$, then $M_H(v)$ is a smooth scheme
of dimension $2 \rho(\br,g,d) = 2(g - r(r-\alpha) )$. If $r \ge 2$ and
$\rho(\br,g,d) < 0$, then $M_H(v)$ is empty. \footnote{Yoshioka proved
that if $\rho \ge 0$, then $M(v)$ is not empty and is irreducible. We
are not using his result.}
\end{thing}
%
%

%
% compare with  $g-k(k-\chi)$, which is a codimension of the classical
% Brill-Noether loci.
%

\begin{thing} {{\bf Lemma (case $r=1$).}}%
\label{r=1}%
  There is an isomorphism
$$
    i_1: \Hilb^d(S) \isom M_H(1,a,\beta),
$$
which takes a 0-dimensional subscheme $\xi \subset S$ of length $d$ to
the isomorphism class $[J_{\xi}(L_a)]$, where $L_a$ is the unique
invertible sheaf on $S$ with Chern class $c_1(L_a)=a$.

  It follows that $M_H(1,a,\beta)$ is a smooth irreducible variety of
dimension $2d$ if $d \ge 0$, and is empty if $d < 0$.
\end{thing}

\begin{thing} {\bf Lemma (case $r=0$).}%
\label{r=0}%
   There is an isomorphism
   $$
      i_0: \Pic^{2g-2-d}(|a|\smirr) \isom M_H(0,a,\beta),
   $$
where $\Pic^{2g-2-d}(|a|\smirr)$ is the ``Picard variety of degree
$2g-2-d$ line bundles on curves in the linear system $|a|^s$''.  A
point in $\Pic^{2g-2-d}(|a|\smirr)$ is a pair $(C,B)$, where $C$ is a
smooth irreducible curve in $|a|$ and $B \in \Pic^{2g-2-d}(C)$. 
In the notations of \cite{grothendieck}, %
$\Pic^{2g-2-d}(|a|\smirr)$ is a component of the
relative Picard variety %
$     %
   {  %
    \mathop{                      %
             \rm \underline{Pic}  %
            } %
    \nolimits %
    }_{ %
            \; {\cal C} / |a|\smirr  %
      } %
$, %
where ${\cal C} / |a|\smirr$ is the universal curve over $|a|\smirr$.
The isomorphism $i_0$ takes a pair $(C,B)$ to the direct image
$(i_C)_*(B)$, where $i_C: C \to B$ is the canonical embedding.

    In particular, $M(0,a,\beta)$ is an irreducible variety of
dimention $2g$. 
%
%(The results of Saint-Donat (\cite{saint-donat}) imply
%that this is not empty even if char k > 0.)
\end{thing}

\begin{thing} {\bf Remark.} The formula $\dim M(r,a,\beta) = 2
\rho(\br,g,d)$ remains valid for $r=0$ and $r=1$.
\end{thing}

\subsection{Brill-Noether stratification of $M_H(v)$ and statement of the %
            main theorem}

\begin{thing}{\bf Lemma:}  %
\label{vanishing-of-h^2}
For every $[F] \in M_H(v)$ we have $H^2(S,F) = 0$.  (Recall that
$(v_1,H) > 0$.)
\end{thing}

This lemma prompts the following definition:
\begin{thing} {\bf Definition:}
$$ BN_k(v) = \{[F]: \ch(F) = v, F \text{ is acceptable and
   $H$-stable}, h^0(S,F) \ge k \}
$$ 
\end{thing}

\begin{thing} {\bf Remark (case $r = 0$).}
If $r = 0$, we have
$$ 
   BN_k(0,a,\beta) =
   W_{2g-2-d}^{\bar{k}}(|a|^s),
$$ where $W_d^r(|a|^s)$ is the relative Brill-Noether scheme over the
linear system $|a|^s$. Note that there is a morphism $W_d^r(|a|^s) \to
|a|^s$ with fibers $W_d^r(C)$.
\end{thing}

% ---------------------------------------------------
%
\begin{thing}{\bf Special 0-subschemes on a surface.}
% Remark (case $r = 1$).}
%
  Let $S$ be an algebraic surface, $L$ be an invertible sheaf on $S$,
and $\xi \in \hilb^d(S)$ be 0-subscheme in $S$. The index of
speciality of $\xi$ with respect to $L$, denoted as $\delta(\xi,L)$,
is defined by the equality
$$
   h^0(S,\jk(L)) =  h^0(S,L) - d + \delta(\xi,L)
$$
It follows that $\delta(\xi,L) \ge 0$, and $\delta(\xi,L) > 0$ if and
only if $\xi$ fails to impose $d$ independent conditions on curves in
the linear system $|L|$. If $\delta(\xi,L) > 0$, $\xi$ is said to be
special with respect to $L$. Considering the adjunction sequence
$$
   \triple{\jk(L)}{L}{\ok(L)},
$$
one can see that $\delta(\xi,L) = h^1(S,\jk(L)) - h^1(S,L)$. 

  For a given $\delta \ge 0$ the set of 0-subschemes on $S$ which
satisfy $\delta(\xi,L) \ge \delta$ form a closed subscheme in the
Hilbert scheme of points $\hilb^d(S)$ which we denote as
$\hilb^d_{(L,\delta)}(S)$. 

\end{thing}
\begin{thing}{\bf Special 0-subschemes on a K3 surface.}
If $S$ is a K3 surface and $L$ is ample, we have $\delta(\xi,L) =
h^1(\jk(L)) = h^0(S,\jk(L)) - \chi(S,\jk(L)) = h^0(S,\jk(L)) + \alpha
- 2$.

{\bf Corollary: case $r = 1$.}
On a K3 surface $S$ there is an equality
$$ BN_k(1,a,\beta) = \hilb^d_{(L,k+\alpha-2)}(S).
$$
\end{thing}

\begin{thing}{\bf Simple Caley-Bacharash 0-cycles.}

{\bf Definition.} A 0-subscheme $\xi \subset S$ of length $d$ is
called simple if it can be written as a sum of $d$ distinct points on
$S$, $\xi = [p_1]+\dots + [p_d]$, $p_i \ne p_j$.

{\bf Definition.} If $L$ is an invertible sheaf on $S$ and $\xi$ is a
simple 0-cycle on $S$, then $\xi$ is said to be Caley-Bacharash with
respect to $L$ if for every point $p_i \in \xi$ and every curve $C$ in
the linear system $|L|$ such that all the points $p_i$, $i \ne d$,
belong to $C$ we also have $p_i \in C$, (i.e., if
$h^0(S,J_{\xi-p_i}(L)) = h^0(S,J_{\xi}(L))$ for every $i$.)

   For example, if $C$ and $D$ are plane curves and $\xi = C \cap D$
is a simple 0-cycle, then $\xi$ is Caley-Bacharash with respect to the
linear system $|\OO_{\PP^2}(\deg C + \deg D - 3)|$. In particular,
every plane cubic containing eight out of nine points of intersection
of two fixed plane cubics also contains the ninth point
(Caley-Bacharash theorem.)

Let $\xi = p_1 + \dots + p_d$ be a simple Caley-Bacharash 0-cycle on
$S$. One can see that $\delta(\xi,L) = \delta(\xi - p_i,L)+ 1$ for
each $i$. It follows that $\xi$ is special with respect to $L$. One
can prove that for each given $\delta > 0$ simple Caley-Bacharash
0-cycles form an open subscheme in $\hilb^d_{(L,\delta)}(S)$.

We let $\hilb^d_{(L,\delta)^o}(S) = \hilb^d_{(L,\delta)}(S) -
\hilb^{d}_{(L,\delta+1)}(S)$.
\end{thing}

%-------------------------------------------

\begin{thing}{\bf Lemma.}
  Let us fix $v$ , and assume that $M(v)$ is not empty. Let us
substitute $M(v)$ with any of its irreducible components \footnote{We
do not use the strong result that $M(v)$ is irreducible}.

\begin{enumerate}
    \item One can construct vector bundles $E$ and $F$ on $M(v)$ and a
map $\sigma: E \to F$ in such a way that $BN_k(v)$ is the set of
points $x \in M(v)$ at which the rank of $\sigma$ drops by a certain
number $l$. In particular, $BN_k(v)$ is the set of points of a closed
subscheme $D_l(\sigma)$ in the moduli space $M(v)$. (Abusing
notations, we sometimes denote this subscheme as $BN_k(v)$). % or
% $BN_k(M(v))$.
%
    \item If $BN_k(v)$ is not empty, then
$$
   \codim_{M(v)} BN_k(M(v)) \le k(k-\chi(v)),
$$
 where $\chi(v) = \chi(S,F)$ for $[F] \in M(v)$. %
\footnote{  %
    Note that the same formula is true in the classical
Brill-Noether theory for line bundles on curves, see, e.g.,
\cite{griffiths-harris}.
         }
\end{enumerate}
\end{thing}

%(Here $\codim_X Y$ for a subscheme $Y$ in an irreducible scheme $X$ is
%the maximum of codimensions of the irreducible components of $Y$).

\begin{thing} {\bf Definition.} We define the expected, or virtual,
dimension of $BN_k(M(v))$ by
\begin{align}
   \vdim BN_k(M(v)) &= \dim M(v) - k(k - \chi(v))     \\
                    &= 2 \rho(\br,g,d) - k(k - \chi)
   \label{eq:virtual-dimension}
\end{align}
where $\chi = \chi(v) = 2r + \beta$ is the Euler characteristics of
any $F \in M_H(v)$.
\end{thing}

Let us fix $S$, $H$, $a \in \Num S$ and $\beta \in \ZZ$. Let $r_0$ be
the positive root of the equation
$$
   \rho(\br_0,g,d) = g - r_0(r_0+\beta) = 0
$$ Note that the moduli space $M_H(r,a,\beta)$ is empty if $r >
r_0$. Let $k_0 = 2r_0 + \beta$ and $k_1 = r_0 + \beta$.

\begin{figure}[htbp]
\begin{center}
\begin{picture}(0,0)%
\includegraphics{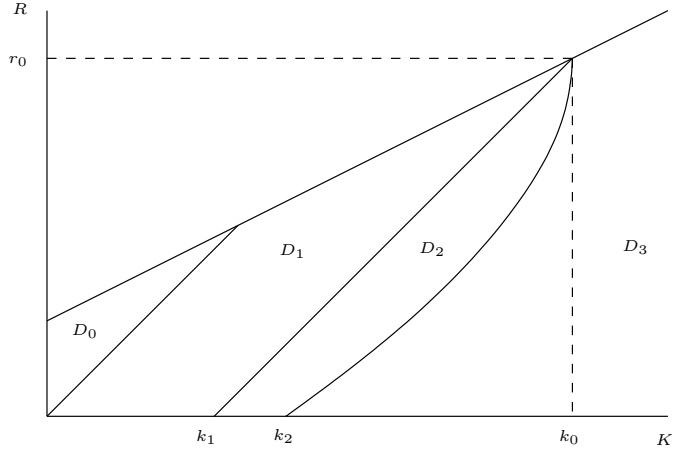}%
% n1
\end{picture}%
\setlength{\unitlength}{3947sp}%
\begingroup\makeatletter\ifx\SetFigFont\undefined%
\gdef\SetFigFont#1#2#3#4#5{%
  \reset@font\fontsize{#1}{#2pt}%
  \fontfamily{#3}\fontseries{#4}\fontshape{#5}%
  \selectfont}%
\fi\endgroup%
\begin{picture}(4153,2856)(360,-2381)
\put(4426,-2341){\makebox(0,0)[lb]{\smash{\SetFigFont{6}{7.2}{\rmdefault}{\mddefault}{\updefault}{\color[rgb]{0,0,0}$K$}%
}}}
\put(384,367){\makebox(0,0)[lb]{\smash{\SetFigFont{6}{7.2}{\rmdefault}{\mddefault}{\updefault}{\color[rgb]{0,0,0}$R$}%
}}}
\put(2941,-1134){\makebox(0,0)[lb]{\smash{\SetFigFont{6}{7.2}{\rmdefault}{\mddefault}{\updefault}{\color[rgb]{0,0,0}$D_2$}%
}}}
\put(4216,-1119){\makebox(0,0)[lb]{\smash{\SetFigFont{6}{7.2}{\rmdefault}{\mddefault}{\updefault}{\color[rgb]{0,0,0}$D_3$}%
}}}
\put(1539,-2318){\makebox(0,0)[lb]{\smash{\SetFigFont{6}{7.2}{\rmdefault}{\mddefault}{\updefault}{\color[rgb]{0,0,0}$k_1$}%
}}}
\put(2025,-2312){\makebox(0,0)[lb]{\smash{\SetFigFont{6}{7.2}{\rmdefault}{\mddefault}{\updefault}{\color[rgb]{0,0,0}$k_2$}%
}}}
\put(360, 59){\makebox(0,0)[lb]{\smash{\SetFigFont{6}{7.2}{\rmdefault}{\mddefault}{\updefault}{\color[rgb]{0,0,0}$r_0$}%
}}}
\put(3819,-2316){\makebox(0,0)[lb]{\smash{\SetFigFont{6}{7.2}{\rmdefault}{\mddefault}{\updefault}{\color[rgb]{0,0,0}$k_0$}%
}}}
\put(2063,-1148){\makebox(0,0)[lb]{\smash{\SetFigFont{6}{7.2}{\rmdefault}{\mddefault}{\updefault}{\color[rgb]{0,0,0}$D_1$}%
}}}
\put(758,-1651){\makebox(0,0)[lb]{\smash{\SetFigFont{6}{7.2}{\rmdefault}{\mddefault}{\updefault}{\color[rgb]{0,0,0}$D_0$}%
}}}
\end{picture}

   \caption{Geography in the $(k,r)$-plane, case $\beta < 0$.}
   \label{fig:geography}
\end{center}
\end{figure}

We define the following domains in the plane with coordinates $(k,r)$:
\begin{align*}
 D_0 &= \{ (k,r): k \ge 0, k \ge 2r + \beta, k < r  \}, \\
 D_1 &= \{ (k,r): r \ge 0, k \ge 2r + \beta, r-k \ge r_0 - k_0  \}, \\
 D_2 &= \{ (k,r): r \ge 0, r-k < r_0 - k_0, \vdim BN_k(r,a,\beta) \ge 0  \}, \\
 D_3 &= \{ (k,r): r \ge 0, k \ge 2r + \beta, \vdim BN_k(r,a,\beta) < 0  \},
\end{align*}
as one the Figure ~\ref{fig:geography}. 

Note that the domain $D_0$ is empty if $\beta < 0$, and
$$
   D_0 \cup D_1 \cup D_2 =
                 \{
                    (k,r):  k \ge \chi(v) = 2r + \beta,
                            \vdim BN_k(r,a,\beta) \ge 0 
                 \}
$$

%
% Consider the following domains in the plane with coordinates $k$ and
%$r$ given by conditions
%\begin{enumerate}
%    \item $k \ge 2r+ \beta$,
%    \item $\rho(\br,g,d) \ge 0$,
%    \item $2 \rho(\br,g,d) - k(k - \chi) \ge 0$.
%\end{enumerate}
%
%The geometric meaning of these conditions is
%\begin{enumerate}
%    \item $h^0(S,F) \ge \chi(S,F)$,
%    \item $\dim M(v) \ge 0$,
%    \item $\vdim BN_k(v) \ge 0$
%\end{enumerate}
%

%

Let 
$$
    V_d^{\br}(|a|^s) = 
                 W_d^{\br}(|a|^s) -
                 W_d^{\overline{r+1}}(|a|^s)
$$
be the subscheme in the relative Picard variety $\Pic^d(|a|^s)$
parametrizing pairs $(C,B)$ for which $h^0(C,B) = r+1$.

Let $(V_d^{\br})'(|a|^s)$ be the open subscheme in $V_d^{\br}(|a|^s)$
parametrizing pairs $(C,B)$ for which $K_C B^{-1}$ is globally
generated, and let $(V_d^{\br})''(|a|^s)$ be the open subscheme in
$V_d^{\br}(|a|^s)$ parametrizing pairs $(C,B)$ for which both $B$ and
$K_C B^{-1}$ are globally generated. 

It is known (the ``existence theorem'' of the classical Brill-Noether
theory; Kempf, Kleiman-Laksov, Fulton-Lazarsfeld) that for a smooth
genus $g$ curve $C$ the scheme $V_d^r(C)$ is not empty whenever
$\rho(\br,d,g) \ge 0$.  However, there are situations when
$\rho(\br,d,g) \ge 0$, but $(V_d^r(C))''$ is empty; the minimal genus
for which this occurs is $g = 5$. For example, if $C$ is a smooth
plane quintic, the variety $V_5^1(C)$ is two-dimensional, but
$(V_5^1)''(C)$ is empty; some other examples are studied in the
section ~\ref{subsection:examples}.

Let us fix a polarized genus $g$ K3 surface $(S,a)$ and integers $r
\ge 0$ and $d \ge 0$ such that $\rho(\br,g,d) \ge 0$. Consider the
following three conditions:
\begin{itemize}
    \item[{\bf $C_1$:}] The scheme $(V_d^{\br})''(|a|^s)$ is not
          empty;
    \item[{\bf $C_2$:}] The linear system $|a|^s$ contains a
          Brill-Noether general curve;
    \item[{\bf $C_3$:}] $\Pic S = \ZZ a$
\end{itemize}

Note that the condition $C_3$ and results of Lazarsfeld imply the
condition $C_2$, and the condition $C_2$ implies the condition $C_1$.

\begin{thing} {\bf Main Theorem.} Let $v = (r,a,\beta) \in K'(S)$ and
$\ep(v) = (r,g,d)$, and assume that $M(r,a,\beta)$ is not
empty. Assume further that either $\Pic(S) = \ZZ a$, or
$(a,H)=1$. Then
\begin{itemize}
   \item[0.]
   If $(k,r) \in D_0$, then either $BN_k(r,a,\beta)$ is empty, or
$$
   \dim BN_k(r,a,\beta) = \vdim BN_k(r,a,\beta);
$$
   \item[1.]
   If $(k,r) \in D_1$, and at least one of the conditions $C_1$, $C_2$
and $C_3$ is satisfied, then $BN_k(r,a,\beta)$ is not empty, and
$$
   \dim BN_k(r,a,\beta) = \vdim BN_k(r,a,\beta);
$$
   \item[2.]
   If $(k,r) \in D_2$ and $r \ge 2$, then $BN_k(r,a,\beta)$ is empty,
even though \\ $\vdim BN_k(r,a,\beta) > 0$;
   \item[3.] 
   If $(k,r) \in D_3$, then $BN_k(r,a,\beta)$ is empty, as expected.
\end{itemize}
\end{thing}

\begin{thing} {\bf Conjecture:}  If $\Pic S = \ZZ a$, then for any
$(k,r) \in D_0$ the variety $BN_k(v)$ is not empty.
\end{thing}

\begin{thing} {\bf Remark.} The second statement could not be
made stronger: we give examples of situations in which $r = 0$ or $r =
1$, $(k,r) \in D_2$, and $BN_k(r,a,\beta)$ is not empty.
\end{thing}

%From now on we sometimes refer to $D_1$ as to a ``good geographic
%region''.

\begin{thing} {\bf Remark.}
Let $k_2$ be the positive root of the equation
$$
   2g - k_2(k_2 - \beta) = 0
$$ The point $(k_2,0)$ is on the boundary of $D_2$, as on the Figure
~\ref{fig:geography}. 
%
% We want to give a few comments upon the
% inequality $k_1 \ne k_2$ and the geography for $r = 0$. 
Note that $k_1$ satisfies the equality
$$
     g - k_1(k_1 - \beta) = 0,
$$
i.e., for a genus $g$ curve $C$ and an integer $k = r+1$ one has $
\vdim W_d^{r}(C) = g - k (k-\beta)$, and 
$$ k \le k_1 \text{ iff }
   \vdim W_d^{r}(C) \ge 0 
$$

Let $\pi: C \to T$ be a family of smooth genus $g$ curves, and let
$\Pic^d_{C/T}$ be the relative Picard variety. One can define the
relative Brill-Noether scheme $(W_d^r)_{C/T}$ over $T$ in such a way
that the fiber of $(W_d^r)_{C/T}$ over a point $t \in T$ is isomorphic
to the Brill-Noether subscheme $W_d^r(C_t)$ in the fiber
$(\Pic^d_{C/T})_t \simeq \Pic^d(C_t)$. (Note that we do not expect
$(W_d^r)_{C/T}$ to be flat over $T$.)  It is natural to define
$$ \vdim (W_d^r)_{C/T} = \dim T + \vdim W_d^r(C_t) = \dim T +
   \rho(r,g,d)
$$

One would expect that if the family $C/T$ is generic enough, then
$$
   \dim (W_d^r)_{C/T} = \vdim (W_d^r)_{C/T}
$$. In particular, if $C/T$ is generic enough and $\dim T = g$, as is
the case when $T = |a|^s$, then one would expect
$$% 
  \dim (W_d^r)_{C/T} = \vdim (W_d^r)_{C/T} = g + (g - k (k-\beta)) =  %
                      2g - k (k-\beta)
$$ It follows that $k_2$ can be defined by the property 
$$ k = r+1 \le k_2 \text{ iff } \vdim (W_d^r)_{C/|a|^s} \ge 0, $$
where $C/|a|^s$ is the universal curve over the linear system $|a|^s$.
In particular, one would expect that for each $k$ such that $k_1 < k =
r+1 < k_2$ the variety $W_d^r(C_t)$ is empty for a generic point $t
\in T$ and that there is a codimension $|\rho(r,g,d)|$ subscheme in
$T$ parametrizing points $t \in T$ such that $W_d^r(C_t)$ is not
empty.  In particular, one would expect that for $k_1 < k = r+1 < k_2$
the variety $(W_d^r)_{C/|a|^s}$ is not empty.

However, if $S$ is a K3 surface and $a \in \Num S$, one can not expect
that the variety $(W_d^r)_{C/|a|^s}$ has the expected dimension, i.e.,
that the classifying map from $|a|^s$ to the moduli space of genus $g$
curves is of general position with respect to the stratification of
$M_g$ by Brill-Noether special curves. (For example, we know that if
one curve in the linear system $|C|$ is hyperelliptic, than all the
curves in $|C|$ are hyperelliptic, and so on.) This discrepancy
creates the region $D_2$ where the Brill-Noether loci are empty in the
regular examples, even though their expected dimension is positive.

\end{thing}

%  is it better to say that ($\ep(v)$ satisfies the condition 
%  ``globally generated and dual is globally generated'')?
%  for every non-hyperelliptic curve $C$ of genus $g$, if $V_d^{\br}(C)$ 
%   is not empty, then $V_d^{\br}(C) - (A_d^{\br}(C) + B_d^{\br}(C))$
% is not empty

\begin{comment}
  The idea behind the proof is to generalize the techniques of
Lazarsfeld (\cite{lazarsfeld-86}, etc.) and Tyurin (~\cite{tyurin}) to
give a set system of correspondences between moduli spaces of
acceptable sheaves with same $c_1$ and $c_2$ but varying rank which
relate different Brill-Noether loci.
\end{comment}

%-------------------------------------------------------------------
%
\subsection{The correspondence $\cA_{r,r'}$}
 
   Let us fix $v \in K'(S)$, $v = (r,v_1,v_2)$, and an integer $r'$
such that $0 \le r' < r$. Let $v' = (r',v_1,v_2)$ and $l = r-r'$.

\begin{thing}{\bf Proposition.}%

 There is a coarse moduli space $A_l(v)$ classifying pairs $(E,V)$,
where $E$ is a coherent sheaf on $S$ and $V \subset H^0(S,E)$ is a
vector space of dimension $l$, which satisfy the following conditions:
\begin{enumerate}
   \item $E$ is acceptable and $H$-stable, 
   \item $\ch(E)=v$,
   \item the canonical ``evaluation'' map 
             $$e_{E,V}: V \tensor \OO_S \to E$$
         is monomorphic,
   \item $\coker e_{E,V}$ is acceptable and $H$-stable.
\end{enumerate}
\end{thing}

The proof of the existence of such a moduli space draws on the results
of Le Potier and He Min (\cite{le-potier,he-min,he-min-1,he-min-2}).

Given a scheme $X$ and morphisms $f: X \to Y$ and $g: X \to Z$, we say
that we are given correspondence between $Y$ and $Z$, and denote this
correspondence as $(X,f,g)$.

The definition of the moduli space $A_l(v)$ gives the correspondence
$\cA_{r,r'}(v) := (A_l(v),\pi_1,\pi_2)$ between $M(v)$ and $M(v')$, 
\begin{equation}
   \correspondence{A_l(v)}{M(v)}{M(v')}{\pi_1}{\pi_2}
   \label{diagram:correspondence}
\end{equation}
where $\pi_1([(E,V)])= [E]$ and $\pi_2([(E,V)])= \coker e_{E,V}$.

A point $[(E,V)] \in A_l(v)$ gives an exact sequence
$$
    0 \to   V \tensor \os \strelka{e_{(E,V)}} E \to F \to 0
$$
Since $h^1(S,O_S)=0$, we have $h^0(S,F) = h^0(S,E) - l$.

Let $A_l^k(v)$ be the closed subscheme in $A_l(v)$ given by the
condition $h^0(S,E) \ge k$.  Let $(A_l^k)^o(v) = A_l^k(v) -
A_l^{k+1}(v)$, and let $BN_k^o(v):= BN_k(v) -
BN_{k+1}(v)$. Restricting the correspondence $\cA_{r,r'}(v)$ to the
locally closed subscheme $(A_l^k)^o(v)$ of $A_l(v)$, we get a
correspondence
$
   \cA_{r,r'}^{k,k'}(v):= ( (A_l^k)^o(v),\pi_1^k,\pi_2^k  ):
$
$$
   \correspondence{(A_l^k)^o(v)}{BN_k^o(v)}{BN_{k-l}^o(v')}{\pi_1^k}{\pi_2^k}
$$

\begin{thing} {\bf Lemma.} For every $[E] \in BN_k(v)^o(v)$ the fiber
$(\pi_1)^{-1}([E])$ is isomorphic to an open subscheme in the
Grassmanian variety $\Gr(l,H^0(S,E))$, and for every $[F] \in
BN_{k-l}(v)^o(v)$ the fiber $(\pi_2)^{-1}([F])$ is isomorphic to an
open subscheme in the Grassmanian variety $\Gr(l,\Ext^1_S(F,\os))$.
\end{thing}

% For a proof, see ~\cite{bn-K3}

%{\bf Remark.} So far we have only used the property $h^1(S,\os)=0$.

\begin{thing}{\bf Remark.} Let $C$ be a smooth algebraic curve, $\Pic^d(C)$
be the Picard variety of degree $d$ line bundles on $C$, and let
$A_1(d)$ be the moduli space of pairs $(L,V)$, where $L \in \Pic^d(C)$
and $V$ is an one-dimensional vector subspace in $H^0(C,L)$. Existence
of such a moduli space follows from the results of Le Potier, or it
can be constructed directly. A point $(L,V) \in A_1(d)$ gives an
exact sequence
$$
    0 \to V \tensor \oc \strelka{e_{(L,V)}} L \to \coker e_{(L,V)} \to 0,
$$
which is isomorphic to the ``adjunction'' sequence
$$
         0 \to \oc \to \oc(D) \to \OO_D(D) \to 0,
$$
where $D$ is the divisor of zeroes of a non-zero section $s \in V$.
Since any effective divisor $D$ on $C$ can be written as a linear
combination of points, $D = \sum n_i [p_i]$, the diagram
\begin{equation}
   \correspondence{A_1(d)(C)}{Pic^d(C)}{C^{(d)}}{\pi_1}{\pi_2}
   \label{diagram:correspondence-curves}
\end{equation}
where $C^{(d)}$ is a symmetric power of the curve $C$, $\pi_1(L,V) =
L$ and $\pi_2(L,V) = \coker e_{L,V}$, can be considered as an analogue
of the correspondence ~\ref{diagram:correspondence}. One can prove
that $\pi_2$ is an isomorphism, and therefore the correspondence
~\ref{diagram:correspondence-curves}%
% $(A_1(d)(C),\pi_1,\pi_2)$
``reduces'' to the Abel-Jacobi morphism $\pi_1\pi_2^{-1}: C^{(d)} \to
Pic^d(C)$.  This example lets us think about two particular
correspondences, $\cA_{r,1}: \;\; \Hilb^d(S) \longmapsto M(r,a,\beta)$
and $\cA_{r,0}: \;\; \Pic^{2g-2-d}(|a|) \longmapsto M(r,a,\beta)$, as
two different analogs of the Abel-Jacobi correspondence in dimension
two, corresponding to the fact that there are two types of subschemes
on a surface: 0-cycles and curves.
\end{thing}

%------------------------------------------------------------
%
\subsection{A few technical lemmas}
\label{technical}

Our study of the correspondence $\cA_{r,r'}$ is based on the following
six technical lemmas:

\begin{thing} {\bf Lemma: Generic extensions are acceptable.}    %1
\label{generic-extensions-are-acceptable}
      Let $S$ be an algebraic surface, $F$ be an acceptable coherent
   sheaf on $S$, and let $W$ be a $k$-vector space. Let $\mathbb A$ be
   the affine scheme associated with the $k$-vector space $\Ext^1_S(F,W
   \tensor \OO_S)$. The affine space $\mathbb A$ parametrizes extension
   classes of the form
   \begin{equation*}
      \triple{W \tensor \os}{E}{F}
   \end{equation*}
   We say that such an extension acceptable if $E$ is acceptable.
   \begin{enumerate}
      \item If $\rk F \ge 2$, i.e., $F$ is locally free, then all such
   extensions are acceptable;
      \item if $\rk F = 1$, i.e., $F = \jk(L)$, where $L$ is an ample
   invertible sheaf on $S$, $\xi$ is a simple effective 0-subscheme on
   $S$, and the pair $(\xi,L+K_S)$ is Caley-Bacharash, then there is a
   {\it nonempty} open subscheme $U \subset \mathbb A$ such that all the
   extensions parametrized by $U$ are acceptable;
      \item if $\rk F = 0$, i.e., $F=(i_C)_* (B)$, then
      \begin{enumerate}%
      \item If $\dim W = 1$, then all the extensions parametrized by
         $\mathbb A - \{0\}$ are acceptable;
      \item If $\dim W \ge 2$, and $A := N_{C/X}B^{-1}$ is a globally
         generated line bundle on $C$, then there is a {\it nonempty} open 
         subscheme $U \subset \mathbb A$ such that every extension class 
         $e \in U$ is acceptable;
      \item In particular, let $\dim W \ge 2$, and assume that 
         $A = N_{C/X}B^{-1}$
         is globally generated and that $\dim W = H^0(C,A)$. 
         Let $\alpha_e: W\dual \to H^0(C,A)$ be the image of $e$ under the 
         isomorphism
         $$
            \Ext^1_S(i_*B,W \tensor \os) \simeq 
            H^1(C,W \tensor A)           \simeq 
            \Hom(W\dual,H^0(C,A))
         $$
         If $\alpha_e$ is an isomorphism, then $e$ is acceptable.
      \end{enumerate}
   \end{enumerate}
%

% Lemma 1.18.

\end{thing}

%-------------------------

\begin{thing} {\bf Lemma : Generic evaluation maps in the globally generated
case are monomorphic.} %2
\label{generic-evaluation-maps-in-globally-generated-case-are-monomorphic}
 Let $X$ be a reduced irreducible scheme, $E$ be a globally generated
locally free sheaf on $X$, and $l \le \rk E$. Then there is a {\it
nonempty} open subset in the Grassmanian variety $Gr(l,H^0(E))$
parametrizing vector subspaces $V \subset H^0(X,E)$ for which the
evaluation map $e_{(E,V)}: V \tensor \ox \to E$ is monomorphic.
%

% , Lemma 2.1.1.

\end{thing}

%-------------------------

\begin{thing} {\bf Lemma: Generic factors in the globally generated case are
                          acceptable.} %3
\label{generic-factors-in-globally-generated-case-are-acceptable}%
   Let $S$ be a smooth algebraic surface, ${\rm char} k = 0$, let $E$ be
   a globally generated locally free sheaf on $S$, and let $1 \le l \le
   \rk E$.  Then there is a {\it nonempty} open subset in the Grassmanian
   variety $Gr(l,H^0(E))$ parametrizing vector subspaces $V \subset
   H^0(S,E)$ for which the cokernel of the evaluation map
   $$e_V: H^0(S,E) \tensor_{k} \OO_S \to E$$ is acceptable.
%

% Proposition 5.

\end{thing}

%-------------------------

\begin{thing}{\bf Lemma: Stability of extensions.} %4
\label{stability-of-extensions} 
   Let $(S,H)$ be a polarized algebraic surface, and let $F$ be an
   $H$-stable acceptable sheaf on $S$. Assume that we are given an
   extension
   \begin{equation*}
      0 \to W \tensor \OO_X \to E \to F \to 0
   \end{equation*}
   Let $e$ be the class of this extension in $\Ext^1(F,W \tensor \OO_X)$,
   and let $\alpha_e$ be the image of $e$ under the isomorphism
   $\Ext^1(F,\triv{W}) \isom \Hom(W\dual,\Ext^1(F,\OO))$.

   \begin{enumerate}
     \item If $\alpha_e$ is not injective, then $E$ is not $H$-stable;
     \item If $\alpha_e$ is injective, $F$ is acceptable, and either
     \begin{enumerate}
     \item
         $\Pic S \simeq \ZZ \cdot c_1(E)$, or
     \item
         $(c_1(F),H) =1$,
   \end{enumerate}
      then $E$ is $H$ - stable.
   \end{enumerate}
   %

% part2, Lemma 1.21.1.

\end{thing}

%-------------------------

\begin{thing}{\bf Lemma: Stability of factors.} %5
\label{stability-of-factors}
  Let $(S,H)$ be a polarized algebraic surface, let $E$ be an
$H$-stable coherent sheaf on $S$, and assume that we are given an
exact sequence of the form
\begin{equation*}
   0 \to W \tensor \OO_S \to E \to F \to 0
\end{equation*}
where $F$ is a torsion-free sheaf.  Assume further that either $\Pic S
\simeq \ZZ \cdot H$ and $c_1(E) = H$, or that $(c_1(E),H) =1$. Then
$F$ is $H$-stable.

%  Lemma 1.21.2.

\end{thing}
%-------------------------

\begin{thing}{\bf Lemma: The correspondence preserves the property of being
globally generated.} %6
\label{correspondence-preserves-globally-generated}
     Let $X$ be a scheme satisfying $h^0(X,\OO_X)=1$, and
   assume that we are given an extension of coherent sheaves
   $$
      \triple{\OO_X}{E}{F}
   $$
    If $E$ is globally generated, then $F$ is
   globally generated. If $F$ is globally generated and $H^1(X,\OO_X)=0$,
   then $E$ is globally generated.
\end{thing}

% , Lemma 1.23.

The proof of lemmas ~\ref{generic-extensions-are-acceptable} -
~\ref{correspondence-preserves-globally-generated} will be given in
\cite{bn-K3-part2}.

%-------------------------------------------------------------------
%
\subsection{Numerical structure of the correspondence $\cA$ on a K3 surface}

Let $S$ be a K3-surface. Consider the map $p$ from the set of the
isomorphism classes of acceptable coherent sheaves on $S$ to the free
Abelian group $P=\ZZ^2$, $p(E) = (h^0(S,E),\rk(E))$. We denote the
coordinates on $P$ as $(K,R)$.

   If $F$ is an acceptable sheaf on $S$, $v=\ch(F)=(r,a,\beta)$ and
$\ep(v)=(r,g,d)$, we have
$$
    \beta = g-1-d
$$
and define
$$
    \alpha = - \beta = d + 1 - g
$$
By the Riemann-Roch formula, we have $\chi(S,F) = 2r+\beta$.

\begin{figure}[!h]
\begin{center}
\setlength{\unitlength}{4144sp}%
\begingroup\makeatletter\ifx\SetFigFont\undefined%
\gdef\SetFigFont#1#2#3#4#5{%
  \reset@font\fontsize{#1}{#2pt}%
  \fontfamily{#3}\fontseries{#4}\fontshape{#5}%
  \selectfont}%
\fi\endgroup%
\begin{picture}(5183,1716)(464,-1345)
{\color[rgb]{0,0,0}\thinlines
\put(1882, 24){\circle{8}}
}%
{\color[rgb]{0,0,0}\put(655,-590){\line( 2, 1){1635.200}}
}%
{\color[rgb]{0,0,0}\put(3393,359){\line( 0,-1){1357}}
\put(3393,-998){\line( 1, 0){1971}}
}%
{\color[rgb]{0,0,0}\put(3858,-998){\line( 2, 1){1498.800}}
}%
{\color[rgb]{0,0,0}\put(655,315){\line( 0,-1){1313}}
\put(655,-998){\line( 1, 0){1976}}
}%
\put(5432,-1157){\makebox(0,0)[lb]{\smash{{\SetFigFont{5}{6.0}{\familydefault}{\mddefault}{\updefault}{\color[rgb]{0,0,0}$K$}%
}}}}
\put(3790,-1135){\makebox(0,0)[lb]{\smash{{\SetFigFont{5}{6.0}{\familydefault}{\mddefault}{\updefault}{\color[rgb]{0,0,0}$\beta$}%
}}}}
\put(2530,-1134){\makebox(0,0)[lb]{\smash{{\SetFigFont{5}{6.0}{\familydefault}{\mddefault}{\updefault}{\color[rgb]{0,0,0}$K$}%
}}}}
\put(485,-590){\makebox(0,0)[lb]{\smash{{\SetFigFont{5}{6.0}{\familydefault}{\mddefault}{\updefault}{\color[rgb]{0,0,0}$\frac{\alpha}{2}$}%
}}}}
\put(1814,-78){\makebox(0,0)[lb]{\smash{{\SetFigFont{5}{6.0}{\familydefault}{\mddefault}{\updefault}{\color[rgb]{0,0,0}$BN_{\chi(v)}^o(v)$}%
}}}}
\put(464,242){\makebox(0,0)[lb]{\smash{{\SetFigFont{5}{6.0}{\familydefault}{\mddefault}{\updefault}{\color[rgb]{0,0,0}$R$}%
}}}}
\put(3173,267){\makebox(0,0)[lb]{\smash{{\SetFigFont{5}{6.0}{\familydefault}{\mddefault}{\updefault}{\color[rgb]{0,0,0}$R$}%
}}}}
\put(1784,229){\makebox(0,0)[lb]{\smash{{\SetFigFont{7}{8.4}{\familydefault}{\mddefault}{\updefault}{\color[rgb]{0,0,0}$l(\beta)$}%
}}}}
\put(1388,-1314){\makebox(0,0)[lb]{\smash{{\SetFigFont{5}{6.0}{\familydefault}{\mddefault}{\updefault}{\color[rgb]{0,0,0}$\beta < 0$}%
}}}}
\put(4204,-1314){\makebox(0,0)[lb]{\smash{{\SetFigFont{5}{6.0}{\familydefault}{\mddefault}{\updefault}{\color[rgb]{0,0,0}$\beta > 0$}%
}}}}
\end{picture}%

   \caption{Geography of the maximal Brill-Noether strata}
   \label{picture:line-l}
\end{center}
\end{figure}
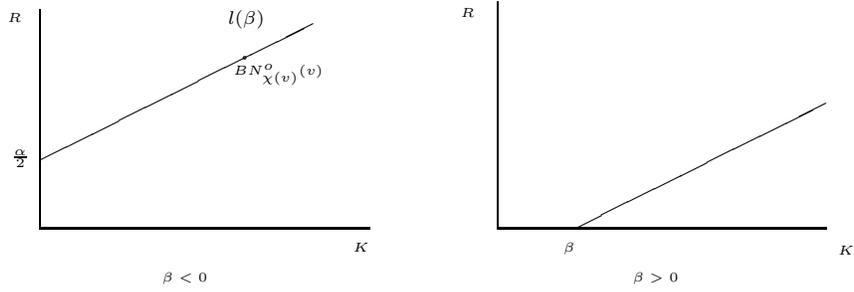

   If $F$ is an acceptable sheaf on $S$ satisfying $h^1(S,F) = 0$ and
$h^2(S,F) = 0$, then $h^0(S,F)=\chi(S,F) = 2r+\beta$, and $p(E)$ is a
point on the line $l(\beta): (K=2R+\beta)$ in $P_{\RR} = P \tensor
\RR$. I.e., $p$ maps the maximal Brill-Noether strata
$BN_{\chi(v)}(v)^o$ to the line $l = l(\beta)$, as on the Figure
~\ref{picture:line-l}. 
%
% (two lines, one for beta < 0, one for beta > 0).
%

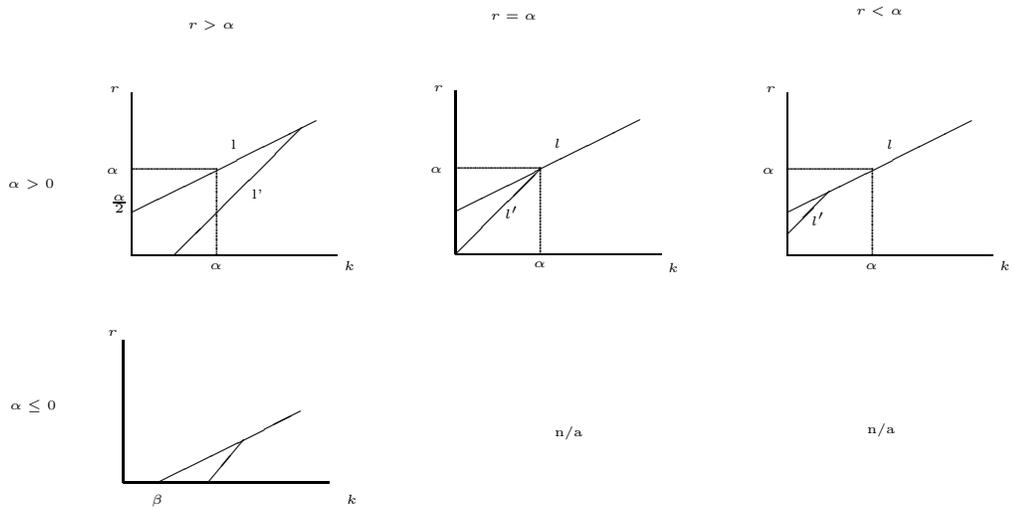
\begin{figure}[!h]
\begin{center}
\setlength{\unitlength}{4144sp}%
\begingroup\makeatletter\ifx\SetFigFont\undefined%
\gdef\SetFigFont#1#2#3#4#5{%
  \reset@font\fontsize{#1}{#2pt}%
  \fontfamily{#3}\fontseries{#4}\fontshape{#5}%
  \selectfont}%
\fi\endgroup%
\begin{picture}(5935,3100)(839,-2486)
\thinlines
{\color[rgb]{0,0,0}\put(1529,-1470){\line( 0,-1){849}}
\put(1529,-2319){\line( 1, 0){1230}}
}%
{\color[rgb]{0,0,0}\put(2250,-2064){\line(-5,-6){212.295}}
}%
{\color[rgb]{0,0,0}\put(1741,-2319){\line( 2, 1){849.600}}
}%
\put(1437,-1440){\makebox(0,0)[lb]{\smash{\SetFigFont{5}{6.0}{\familydefault}{\mddefault}{\updefault}{\color[rgb]{0,0,0}$r$}%
}}}
\put(1698,-2446){\makebox(0,0)[lb]{\smash{\SetFigFont{5}{6.0}{\familydefault}{\mddefault}{\updefault}{\color[rgb]{0,0,0}$\beta$}%
}}}
\put(2867,-2446){\makebox(0,0)[lb]{\smash{\SetFigFont{5}{6.0}{\familydefault}{\mddefault}{\updefault}{\color[rgb]{0,0,0}$k$}%
}}}
{\color[rgb]{0,0,0}\put(1578, 17){\line( 0,-1){977}}
\put(1578,-960){\line( 1, 0){1230}}
}%
{\color[rgb]{0,0,0}\put(1578,-705){\line( 2, 1){1104}}
}%
{\color[rgb]{0,0,0}\put(2597,-195){\line(-1,-1){765}}
}%
{\color[rgb]{0,0,0}\multiput(1578,-451)(8.92982,0.00000){58}{\makebox(1.5875,11.1125){\SetFigFont{5}{6}{\rmdefault}{\mddefault}{\updefault}.}}
\multiput(2087,-451)(0.00000,-8.92982){58}{\makebox(1.5875,11.1125){\SetFigFont{5}{6}{\rmdefault}{\mddefault}{\updefault}.}}
}%
{\color[rgb]{0,0,0}\put(3515, 25){\line( 0,-1){977}}
\put(3515,-952){\line( 1, 0){1231}}
}%
{\color[rgb]{0,0,0}\put(3515,-697){\line( 2, 1){1104}}
}%
{\color[rgb]{0,0,0}\put(3515,-952){\line( 1, 1){509.500}}
}%
{\color[rgb]{0,0,0}\multiput(3515,-443)(8.94737,0.00000){58}{\makebox(1.5875,11.1125){\SetFigFont{5}{6}{\rmdefault}{\mddefault}{\updefault}.}}
\multiput(4025,-443)(0.00000,-8.92982){58}{\makebox(1.5875,11.1125){\SetFigFont{5}{6}{\rmdefault}{\mddefault}{\updefault}.}}
}%
{\color[rgb]{0,0,0}\put(5501, 17){\line( 0,-1){977}}
\put(5501,-960){\line( 1, 0){1230}}
}%
{\color[rgb]{0,0,0}\put(5501,-705){\line( 2, 1){1103.200}}
}%
{\color[rgb]{0,0,0}\put(5755,-578){\line(-1,-1){254}}
}%
{\color[rgb]{0,0,0}\multiput(5501,-451)(8.92982,0.00000){58}{\makebox(1.5875,11.1125){\SetFigFont{5}{6}{\rmdefault}{\mddefault}{\updefault}.}}
\multiput(6010,-451)(0.00000,-8.92982){58}{\makebox(1.5875,11.1125){\SetFigFont{5}{6}{\rmdefault}{\mddefault}{\updefault}.}}
}%
\put(1428,-471){\makebox(0,0)[lb]{\smash{\SetFigFont{5}{6.0}{\familydefault}{\mddefault}{\updefault}{\color[rgb]{0,0,0}$\alpha$}%
}}}
\put(2047,-1041){\makebox(0,0)[lb]{\smash{\SetFigFont{5}{6.0}{\familydefault}{\mddefault}{\updefault}{\color[rgb]{0,0,0}$\alpha$}%
}}}
\put(1450, 17){\makebox(0,0)[lb]{\smash{\SetFigFont{5}{6.0}{\familydefault}{\mddefault}{\updefault}{\color[rgb]{0,0,0}$r$}%
}}}
\put(2851,-1045){\makebox(0,0)[lb]{\smash{\SetFigFont{5}{6.0}{\familydefault}{\mddefault}{\updefault}{\color[rgb]{0,0,0}$k$}%
}}}
\put(1450,-662){\makebox(0,0)[lb]{\smash{\SetFigFont{5}{6.0}{\familydefault}{\mddefault}{\updefault}{\color[rgb]{0,0,0}$\frac{\alpha}{2}$}%
}}}
\put(839,-557){\makebox(0,0)[lb]{\smash{\SetFigFont{5}{6.0}{\familydefault}{\mddefault}{\updefault}{\color[rgb]{0,0,0}$\alpha > 0$}%
}}}
\put(2172,-323){\makebox(0,0)[lb]{\smash{\SetFigFont{5}{6.0}{\rmdefault}{\mddefault}{\updefault}{\color[rgb]{0,0,0}l}%
}}}
\put(2299,-620){\makebox(0,0)[lb]{\smash{\SetFigFont{5}{6.0}{\rmdefault}{\mddefault}{\updefault}{\color[rgb]{0,0,0}l'}%
}}}
\put(1917,399){\makebox(0,0)[lb]{\smash{\SetFigFont{5}{6.0}{\rmdefault}{\mddefault}{\updefault}{\color[rgb]{0,0,0}$r > \alpha$}%
}}}
\put(3366,-463){\makebox(0,0)[lb]{\smash{\SetFigFont{5}{6.0}{\familydefault}{\mddefault}{\updefault}{\color[rgb]{0,0,0}$\alpha$}%
}}}
\put(3985,-1033){\makebox(0,0)[lb]{\smash{\SetFigFont{5}{6.0}{\familydefault}{\mddefault}{\updefault}{\color[rgb]{0,0,0}$\alpha$}%
}}}
\put(3387, 25){\makebox(0,0)[lb]{\smash{\SetFigFont{5}{6.0}{\familydefault}{\mddefault}{\updefault}{\color[rgb]{0,0,0}$r$}%
}}}
\put(4109,-315){\makebox(0,0)[lb]{\smash{\SetFigFont{5}{6.0}{\rmdefault}{\mddefault}{\updefault}{\color[rgb]{0,0,0}$l$}%
}}}
\put(3812,-740){\makebox(0,0)[lb]{\smash{\SetFigFont{5}{6.0}{\rmdefault}{\mddefault}{\updefault}{\color[rgb]{0,0,0}$l'$}%
}}}
\put(3728,450){\makebox(0,0)[lb]{\smash{\SetFigFont{5}{6.0}{\rmdefault}{\mddefault}{\updefault}{\color[rgb]{0,0,0}$r = \alpha$}%
}}}
\put(4788,-1060){\makebox(0,0)[lb]{\smash{\SetFigFont{5}{6.0}{\familydefault}{\mddefault}{\updefault}{\color[rgb]{0,0,0}$k$}%
}}}
\put(5352,-471){\makebox(0,0)[lb]{\smash{\SetFigFont{5}{6.0}{\familydefault}{\mddefault}{\updefault}{\color[rgb]{0,0,0}$\alpha$}%
}}}
\put(5970,-1041){\makebox(0,0)[lb]{\smash{\SetFigFont{5}{6.0}{\familydefault}{\mddefault}{\updefault}{\color[rgb]{0,0,0}$\alpha$}%
}}}
\put(5373, 17){\makebox(0,0)[lb]{\smash{\SetFigFont{5}{6.0}{\familydefault}{\mddefault}{\updefault}{\color[rgb]{0,0,0}$r$}%
}}}
\put(6774,-1045){\makebox(0,0)[lb]{\smash{\SetFigFont{5}{6.0}{\familydefault}{\mddefault}{\updefault}{\color[rgb]{0,0,0}$k$}%
}}}
\put(6095,-323){\makebox(0,0)[lb]{\smash{\SetFigFont{5}{6.0}{\rmdefault}{\mddefault}{\updefault}{\color[rgb]{0,0,0}$l$}%
}}}
\put(5646,-778){\makebox(0,0)[lb]{\smash{\SetFigFont{5}{6.0}{\rmdefault}{\mddefault}{\updefault}{\color[rgb]{0,0,0}$l'$}%
}}}
\put(5978,-2025){\makebox(0,0)[lb]{\smash{\SetFigFont{5}{6.0}{\rmdefault}{\mddefault}{\updefault}{\color[rgb]{0,0,0}n/a}%
}}}
\put(4111,-2043){\makebox(0,0)[lb]{\smash{\SetFigFont{5}{6.0}{\rmdefault}{\mddefault}{\updefault}{\color[rgb]{0,0,0}n/a}%
}}}
\put(851,-1884){\makebox(0,0)[lb]{\smash{\SetFigFont{5}{6.0}{\familydefault}{\mddefault}{\updefault}{\color[rgb]{0,0,0}$\alpha \le 0$}%
}}}
\put(5914,482){\makebox(0,0)[lb]{\smash{\SetFigFont{5}{6.0}{\rmdefault}{\mddefault}{\updefault}{\color[rgb]{0,0,0}$r < \alpha$}%
}}}
\end{picture}

   \caption{Relative position of the lines $l$ and $l'$}
   \label{lines-l-and-l'}
\end{center}
\end{figure}

Note that the correspondences $\cA_{r,r'}^{k,k'}$ acts along the lines
$l'(c): (K-R = c)$ in $P_{\RR}$. Note also that for every $(k',r') \in
P$ the line $l'(k'-r')$ contains $(k',r')$ and intersects the line
$l(\beta)$ at some point $(k,r)$, i.e., for every $(k,r) \in P$ there
is a $(k',r') \in l(\beta)$ and a correspondence $\cA_{r,r'}^{k,k'}$
between $BN_{k'}(r',a,\beta)$ and the maximal Brill-Noether loci
$BN_{k}(r,a,\beta)$, as on the Figure ~\ref{lines-l-and-l'}.

%The following picture demonstrates the four possible combinatorial
%variants for the relative position of $l(\beta)$ and $l'(c)$ in the
%$(K,R)$-plane:

%

Thus the correspondence $\cA$ relates the geometry of Brill-Noether
loci on the lines $l'$. We use this correspondence twice, first time
to derive certain properties of the maximal Brill-Noether loci
$BN_{2r+\beta}(r,a,\beta)$ and second time to relate
$BN_{2r+\beta}(r,a,\beta)$ and a given Brill-Noether variety
$BN_{k'}(r',a,\beta)$ in order to establish that the latter behave in
the way predicted by the formula ~\ref{eq:virtual-dimension} in a
certain domain.

%
%This correspondence allows us to relate various properties of the
%Brill-Noether loci $BN_k(M(r,a,\beta))$ for fixed $a$ and $\beta$ and
%varying $r$ and $k$. Using this correspondence, we establish first the
%nonemptiness of some maximal Brill-Noether loci $BN_{\chi(v)}(M(v))$
%for $r \ge 2$ and the existence of the globally generated invertible
%sheaves in some of these strata. Using this correspondence again, we
%are able to prove the nonemptiness of some of the Brill-Noether loci
%$BN_k(M(r',a,\beta))$, and using the Mukai formula for the dimension
%of $M(v)$ we are able to prove that many of the Brill-Noether loci
%$BN_k(M(v))$ are of the expected dimension.
%
%  In the particular case $r'=0$ we are able to restore the original
%Lazarsfeld's construction (cf \cite{lazarsfeld-86}).
%

%---------------------------------------------------------------
\subsection{Four special cases.}
Let $v = (r,a,\beta)$ and $v' = (r',a,\beta)$, $v, v' \in K'(S)$.  We
will be studying the correspondence $\cA_{r,r'}^{k,k'}$. The following
four cases are particularly interesting:
\begin{itemize}
  \item[(a)] $p(v) \in l$, $r' = 0$; 
  \item[(b)] $p(v) \in l$, $r'$ is any integer such that  $0 \le r' < r$;
  \item[(c)] $r' = 0$;
  \item[(d)] $r' = 1$. 
%  \item[(e)] $r = 1$, $r' = 0$. 
%
\end{itemize}
(Note that (a) is a special case of (b) and (c).)

\begin{figure}[!h]
\begin{center}
\setlength{\unitlength}{4144sp}%
\begingroup\makeatletter\ifx\SetFigFont\undefined%
\gdef\SetFigFont#1#2#3#4#5{%
  \reset@font\fontsize{#1}{#2pt}%
  \fontfamily{#3}\fontseries{#4}\fontshape{#5}%
  \selectfont}%
\fi\endgroup%
\begin{picture}(5252,1267)(1450,-856)
{\color[rgb]{0,0,0}\thinlines
\put(5980,-582){\circle{10}}
}%
{\color[rgb]{0,0,0}\put(6190,-373){\circle{10}}
}%
{\color[rgb]{0,0,0}\put(4514,-653){\circle{10}}
}%
{\color[rgb]{0,0,0}\put(4794,-373){\circle{10}}
}%
{\color[rgb]{0,0,0}\put(3397,-373){\circle{10}}
}%
{\color[rgb]{0,0,0}\put(3677,-94){\circle{10}}
}%
{\color[rgb]{0,0,0}\put(1723,-653){\circle{10}}
}%
{\color[rgb]{0,0,0}\put(2281,-94){\circle{10}}
}%
{\color[rgb]{0,0,0}\put(2979,114){\line( 0,-1){767}}
\put(2979,-653){\line( 1, 0){920}}
}%
{\color[rgb]{0,0,0}\put(2979,-444){\line( 2, 1){837.600}}
}%
{\color[rgb]{0,0,0}\put(3677,-94){\line(-1,-1){279.500}}
}%
{\color[rgb]{0,0,0}\put(1583,114){\line( 0,-1){767}}
\put(1583,-653){\line( 1, 0){921}}
}%
{\color[rgb]{0,0,0}\put(1583,-444){\line( 2, 1){837.600}}
}%
{\color[rgb]{0,0,0}\put(2281,-94){\line(-1,-1){558.500}}
}%
{\color[rgb]{0,0,0}\put(5769,116){\line( 0,-1){767}}
\put(5769,-651){\line( 1, 0){921}}
}%
{\color[rgb]{0,0,0}\put(5769,-442){\line( 2, 1){838}}
}%
{\color[rgb]{0,0,0}\put(6188,-371){\line(-1,-1){210}}
}%
{\color[rgb]{0,0,0}\put(4378,119){\line( 0,-1){767}}
\put(4378,-648){\line( 1, 0){921}}
}%
{\color[rgb]{0,0,0}\put(4378,-439){\line( 2, 1){837.600}}
}%
{\color[rgb]{0,0,0}\put(4797,-369){\line(-1,-1){279}}
}%
\put(2846,106){\makebox(0,0)[lb]{\smash{\SetFigFont{5}{6.0}{\familydefault}{\mddefault}{\updefault}{\color[rgb]{0,0,0}$r$}%
}}}
\put(3407,-157){\makebox(0,0)[lb]{\smash{\SetFigFont{5}{6.0}{\rmdefault}{\mddefault}{\updefault}{\color[rgb]{0,0,0}l}%
}}}
\put(2537,-717){\makebox(0,0)[lb]{\smash{\SetFigFont{5}{6.0}{\familydefault}{\mddefault}{\updefault}{\color[rgb]{0,0,0}$k$}%
}}}
\put(1450,106){\makebox(0,0)[lb]{\smash{\SetFigFont{5}{6.0}{\familydefault}{\mddefault}{\updefault}{\color[rgb]{0,0,0}$r$}%
}}}
\put(2010,-157){\makebox(0,0)[lb]{\smash{\SetFigFont{5}{6.0}{\rmdefault}{\mddefault}{\updefault}{\color[rgb]{0,0,0}l}%
}}}
\put(5637,108){\makebox(0,0)[lb]{\smash{\SetFigFont{5}{6.0}{\familydefault}{\mddefault}{\updefault}{\color[rgb]{0,0,0}$r$}%
}}}
\put(6197,-155){\makebox(0,0)[lb]{\smash{\SetFigFont{5}{6.0}{\rmdefault}{\mddefault}{\updefault}{\color[rgb]{0,0,0}l}%
}}}
\put(6619,-762){\makebox(0,0)[lb]{\smash{\SetFigFont{5}{6.0}{\familydefault}{\mddefault}{\updefault}{\color[rgb]{0,0,0}$k$}%
}}}
\put(4246,109){\makebox(0,0)[lb]{\smash{\SetFigFont{5}{6.0}{\familydefault}{\mddefault}{\updefault}{\color[rgb]{0,0,0}$r$}%
}}}
\put(4806,-155){\makebox(0,0)[lb]{\smash{\SetFigFont{5}{6.0}{\rmdefault}{\mddefault}{\updefault}{\color[rgb]{0,0,0}l}%
}}}
\put(3932,-717){\makebox(0,0)[lb]{\smash{\SetFigFont{5}{6.0}{\familydefault}{\mddefault}{\updefault}{\color[rgb]{0,0,0}$k$}%
}}}
\put(2013,324){\makebox(0,0)[lb]{\smash{\SetFigFont{8}{9.6}{\familydefault}{\mddefault}{\updefault}{\color[rgb]{0,0,0}a}%
}}}
\put(3444,324){\makebox(0,0)[lb]{\smash{\SetFigFont{8}{9.6}{\familydefault}{\mddefault}{\updefault}{\color[rgb]{0,0,0}b}%
}}}
\put(4805,340){\makebox(0,0)[lb]{\smash{\SetFigFont{8}{9.6}{\familydefault}{\mddefault}{\updefault}{\color[rgb]{0,0,0}c}%
}}}
\put(6169,329){\makebox(0,0)[lb]{\smash{\SetFigFont{8}{9.6}{\familydefault}{\mddefault}{\updefault}{\color[rgb]{0,0,0}d}%
}}}
\put(4728,-820){\makebox(0,0)[lb]{\smash{\SetFigFont{5}{6.0}{\familydefault}{\mddefault}{\updefault}{\color[rgb]{0,0,0}$r'=0$}%
}}}
\put(6137,-815){\makebox(0,0)[lb]{\smash{\SetFigFont{5}{6.0}{\familydefault}{\mddefault}{\updefault}{\color[rgb]{0,0,0}$r'=1$}%
}}}
\end{picture}

   \caption{Four special cases}
   \label{picture:4-cases}
\end{center}
\end{figure}
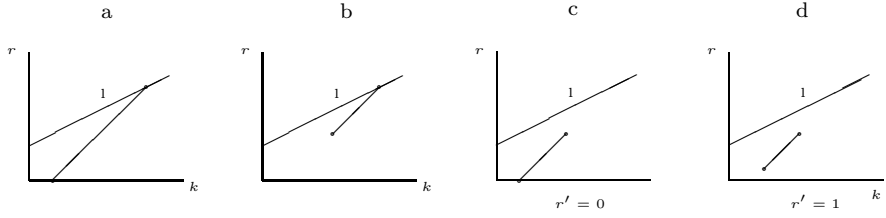

These cases are illustrated by the Figure ~\ref{picture:4-cases}.

We use cases (a) and (b) to prove our main theorem, and cases (c) and
(d) to establish certain birational properties of the birational
geometry of the moduli spaces.

  There is a fifth interesting case, one with $r=1$ and $r'=0$, which
relates the geometry of the Hilbert scheme of points on $S$ and the
relative Picard variety for curves in the linear system $|a|$. We do
not consider it in this text (see ~\cite{bn-K3-part2} for more
details).

 We will analyze these four situations case-by case.
%
%---------------------------------------------------------------------
%
%
\subsubsection{Case (a)}  
\label{section:case-a}
  Let us fix $v=(r,a,\beta)$, $r \ge 2$, and let $v'=(0,a,\beta)$.
Let $k = \chi(v) = 2r+\beta$ and $k' = r+\beta$. Assume that $k \ge
r$, i.e., $(k,r) \notin D_0$. It follows that $k' \ge 0$ and $r \ge
\alpha$.

\begin{figure}[!h]
\begin{center}
\begin{picture}(0,0)%
\includegraphics{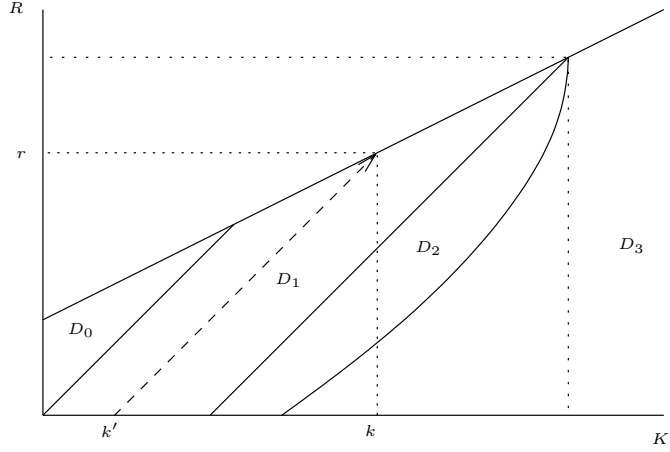}% 
% n2
\end{picture}%
\setlength{\unitlength}{3947sp}%
\begingroup\makeatletter\ifx\SetFigFont\undefined%
\gdef\SetFigFont#1#2#3#4#5{%
  \reset@font\fontsize{#1}{#2pt}%
  \fontfamily{#3}\fontseries{#4}\fontshape{#5}%
  \selectfont}%
\fi\endgroup%
\begin{picture}(4129,2856)(384,-2381)
\put(4426,-2341){\makebox(0,0)[lb]{\smash{\SetFigFont{6}{7.2}{\rmdefault}{\mddefault}{\updefault}{\color[rgb]{0,0,0}$K$}%
}}}
\put(384,367){\makebox(0,0)[lb]{\smash{\SetFigFont{6}{7.2}{\rmdefault}{\mddefault}{\updefault}{\color[rgb]{0,0,0}$R$}%
}}}
\put(2941,-1134){\makebox(0,0)[lb]{\smash{\SetFigFont{6}{7.2}{\rmdefault}{\mddefault}{\updefault}{\color[rgb]{0,0,0}$D_2$}%
}}}
\put(4216,-1119){\makebox(0,0)[lb]{\smash{\SetFigFont{6}{7.2}{\rmdefault}{\mddefault}{\updefault}{\color[rgb]{0,0,0}$D_3$}%
}}}
\put(758,-1651){\makebox(0,0)[lb]{\smash{\SetFigFont{6}{7.2}{\rmdefault}{\mddefault}{\updefault}{\color[rgb]{0,0,0}$D_0$}%
}}}
\put(435,-541){\makebox(0,0)[lb]{\smash{\SetFigFont{6}{7.2}{\rmdefault}{\mddefault}{\updefault}{\color[rgb]{0,0,0}$r$}%
}}}
\put(2625,-2288){\makebox(0,0)[lb]{\smash{\SetFigFont{6}{7.2}{\rmdefault}{\mddefault}{\updefault}{\color[rgb]{0,0,0}$k$}%
}}}
\put(967,-2295){\makebox(0,0)[lb]{\smash{\SetFigFont{6}{7.2}{\rmdefault}{\mddefault}{\updefault}{\color[rgb]{0,0,0}$k'$}%
}}}
\put(2063,-1336){\makebox(0,0)[lb]{\smash{\SetFigFont{6}{7.2}{\rmdefault}{\mddefault}{\updefault}{\color[rgb]{0,0,0}$D_1$}%
}}}
\end{picture}

   \caption{Geography in the case (a).}
   \label{fig:case-a}
\end{center}
\end{figure}

\begin{comment}
  Note that
  $$
     BN_{k'}(0,a,\beta) = 
     BN_{k'}(\Pic^{2g-2-d}(|a|^s)) = 
     V_{2g-2-d}^{k'}(|a|^s)),
  $$
  where
  $$
      V_d^r(|a|^s) =
      W_d^r(|a|^s) -  W_d^{r+1}(|a|^s)
  $$
\end{comment}
Consider the correspondence $\cA_{r,0}^{2r+\beta,r+\beta}$,
$$
   \correspondence{BN_{\chi(v)}^o(A_r(v))}{BN_{\chi(v)}^o(v)}{ %
        V_{2g-2-d}^{\br + \beta}(|a|\smirr)}{\pi_1}{\pi_2}
$$
as illustrated on the Figure ~\ref{fig:case-a} in the case $\beta <
0$.

Note that for a curve $C \in |a|^s$ there is a Serre duality
isomorphism $V_{2g-2-d}^{\br+\beta}(C) \isom V_d^{\br}(C)$.

If $BN_{\chi(v)}^0(v)$ is not empty, then the fiber $\pi_1^{-1}([E])$
over $[E] \in BN_{\chi(v)}^0(v)$ is isomorphic to an open subscheme in
the Grassmanian variety $\Gr(r,H^0(S,E)) \simeq \Gr(r,2r+\beta)$ of
dimension $r(r+\beta)$. (Note that $r+\beta \ge 0$). If
$V_{2g-2-d}^{\br + \beta}(|a|^s)$ is not empty, then the fiber
$\pi_2^{-1}((i_C)_*(B))$ over a point $(C,B) \in
V_{2g-2-d}^{\br + \beta}(|a|^s)$ is isomorphic to an open subscheme in
the Grassmanian variety $\Gr(r,\Ext^1_S(i_*B,\os)) \simeq
\Gr(r,H^0(C,A)) \simeq \Gr(r,r) \simeq \Spec k$, where $A = K_C
B^{-1}$. This proves that $\pi_2$ is an embedding. Note that $h^0(C,A)
= h^1(C,B) = h^0(C,B) - \chi(C,B) = r+\beta-\beta = r$.

\begin{thing} {\bf Lemma.}%
\label{conclusion-of-case-a}%
 Assume that the following conditions are satisfied:
\begin{enumerate}
   \item
       $\rho(\br,g,d) \ge 0$, i.e., $(k,r) \in D_1$,
   \item 
       There is a globally generated line bundle $A \in V_d^{\br}(|a|^s)$, and
   \item 
       Either $\Pic S = \ZZ a$, or $(a,H) = 1$.
\end{enumerate}
Then $BN_{2r+\beta}^o(r,a,\beta)$ is not empty. Moreover, if
$(V_d^{\br})''(|a|^s)$ is not empty, then $BN_{2r+\beta}^o(r,a,\beta)$
contains a globally generated vector bundle.
\end{thing}

{\bf Proof:} By the assumption, there is a $C \in |a|^s$ and a
globally generated line bundle $A \in V_d^{\br}(C)$.  Let $B=K_C
A^{-1}$, and let $W$ be a vector space of dimension $r$.  Since
$$
     \dim_k \Ext^1_S((i_C)_*B,\os) = \dim_k H^0(C,A) = r,
$$
there is a unique (up to an automorphism of $W$) extension class $e$
of the form
$$
      0 \to W \tensor \os \to E \to i_*B  \to 0
$$
such that the induced map $\alpha_e: W \dual \to \Ext^1_S(i_*B,\os)$
is an isomorphism. By Lemma ~\ref{generic-extensions-are-acceptable},
$E$ is locally free. Since either $\Pic S = \ZZ a$ or $(a,H) = 1$,
Lemma ~\ref{stability-of-extensions} implies that $E$ is
$H$-stable. Thus $E \in BN_{\chi(v)}^o(v)$. Moreover, if $B$ is
globally generated, then, by Lemma
~\ref{correspondence-preserves-globally-generated}, $E$ is globally
generated, which proves the lemma.

% -----------------------------
%
\begin{comment}
  \begin{thing}
   Let us study the image of $\pi_2$. Let $A_d^r(|a|^s)$ be a subscheme
  in $V_d^k(|a|^s)$ parametrizing pairs $(C,L)$ such that $L$ is not
  globally generated, and let $B_d^k(|a|^s)$ be a subscheme in
  $V_d^r(|a|^s)$ parametrizing pairs $(C,L)$ such that $K_C L^{-1}$ is
  not globally generated.
  \end{thing}
\end{comment}

Our previous considerations imply the following
\begin{thing}{\bf Lemma:} If $\Pic S = \ZZ a$ or $(a,H) = 1$, then 
  $\pi_2$ is an isomorphism of $BN_{\chi(v)}^o(A_r(v))$ unto
  $(V_{2g-2-d}^{\br + \beta})'(|a|^s)$.
\end{thing}

In particular, the correspondence $\cA_{r,0}^{2r+\beta,r+\beta}$ gives
a morphism
$$
   L: \; \; (V_{2g-2-d}^{\br+\beta})'(|a|^s)
            \to BN_{\chi(v)}^o(v)
$$

\begin{thing} The construction above is equivalent to the
construction given by Lazarsfeld in \cite{lazarsfeld-86}, as we will
now explain.
\end{thing}

    Let $(C,A) \in V_{d}^{\br}(|a|^s)$, and assume that $A$ is
globally generated. The canonical evaluation map $H^0(C,A) \tensor_k
\oc \to A$ induces an epimorphic map $H^0(C,A) \tensor_k \os \to
(i_C)_*A$; denoting its kernel as $F$, we get an exact sequence of
$\os$-modules
$$
   0 \to F \to H^0(C,A) \tensor_k \os \to i_*A \to 0 
$$ 
Note that $F$ is locally free since locally it is a kernel of an
epimorphism $R^n \to R/f \to 0$. 
%, and therefore $\Tor_1(F,k(p))=0$.
By assumptions on $A$ we have $\dim_k H^0(C,A) = r$.  Applying the
functor $\sheafHom_S(\cdot,\os)$ to the exact sequence above, we get
an extension
$$
   0 \to H^0(C,A)\dual \tensor_k \os \to F\dual \to i_*B \to 0
$$
The long corresponding long exact sequence
\begin{gather*}
   0 \to H^0(C,A)\dual \to H^0(S,F\dual) \to H^0(C,B) \to 0, \\
   0 \to H^1(S,F\dual) \to H^1(C,B) \strelka{\delta} H^0(C,A)\dual \to H^2(S,F\dual) \to 0
\end{gather*}
gives
$
  h^0(S,F\dual) = h^0(C,A)+h^0(C,B)= r+(r+\beta) = \chi(v).
$
One can check that $\delta$ is the Serre duality isomorphism, and
therefore $h^1(S,F\dual) = h^2(S,F\dual) = 0$.

If $\Pic S = \ZZ a$ or $(C,H) = 1$, then $F\dual$ is stable by lemma %
~\ref{stability-of-extensions}. %

%   Or, in the opposite logical direction, since $(C,H)=1$, 
% $F\dual$ is stable by lemma (*), 
% and therefore $h^1(S,F\dual) = h^2(S,F\dual) = 0$. 

It follows that $F\dual \in BN_{\chi(v)}^o(r,a,\beta)$.  Thus 
%
%the correspondence $\cA_{r,0}^{2r+\beta,r+\beta}$ ``reduces'' to a
%morphism
we have a morphism
\begin{gather*}
    L: \; \; (V_{2g-2-d}^{\br+\beta})'(|a|^s) \to
    BN_{\chi(v)}^o(r,a,\beta), \\
    (C,B) \mapsto F\dual
\end{gather*}

It follows that the case (a) of the correspondence $\cA$ restores the
construction of Lazarsfeld (\cite{lazarsfeld-86}).

\subsubsection{Case (b) and proof of the main theorem}  
\label{section:case-b}

Let us fix $v= (r,a,\beta)$, $v'=(r',a,\beta)$, $r \ge 2$, % 
$0 \le r' < r$, and $k' \ge \chi(v')$. Let $k = \chi(v) = 2r+\beta$.
Assume that $k \ge r$, i.e., $(k,r) \notin D_0$.

\begin{figure}[htbp]
\begin{center}
\begin{picture}(0,0)%
\includegraphics{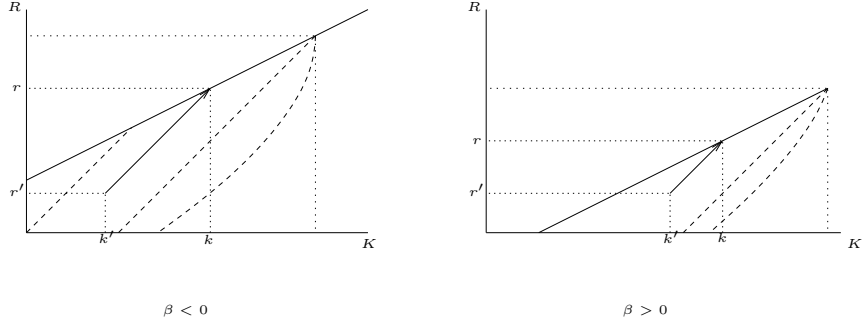}%
% n3
\end{picture}%
\setlength{\unitlength}{2171sp}%
\begingroup\makeatletter\ifx\SetFigFont\undefined%
\gdef\SetFigFont#1#2#3#4#5{%
  \reset@font\fontsize{#1}{#2pt}%
  \fontfamily{#3}\fontseries{#4}\fontshape{#5}%
  \selectfont}%
\fi\endgroup%
\begin{picture}(9585,3606)(384,-3131)
\put(4426,-2341){\makebox(0,0)[lb]{\smash{\SetFigFont{5}{6.0}{\rmdefault}{\mddefault}{\updefault}{\color[rgb]{0,0,0}$K$}%
}}}
\put(384,367){\makebox(0,0)[lb]{\smash{\SetFigFont{5}{6.0}{\rmdefault}{\mddefault}{\updefault}{\color[rgb]{0,0,0}$R$}%
}}}
\put(435,-541){\makebox(0,0)[lb]{\smash{\SetFigFont{5}{6.0}{\rmdefault}{\mddefault}{\updefault}{\color[rgb]{0,0,0}$r$}%
}}}
\put(2625,-2288){\makebox(0,0)[lb]{\smash{\SetFigFont{5}{6.0}{\rmdefault}{\mddefault}{\updefault}{\color[rgb]{0,0,0}$k$}%
}}}
\put(1425,-2280){\makebox(0,0)[lb]{\smash{\SetFigFont{5}{6.0}{\rmdefault}{\mddefault}{\updefault}{\color[rgb]{0,0,0}$k'$}%
}}}
\put(405,-1748){\makebox(0,0)[lb]{\smash{\SetFigFont{5}{6.0}{\rmdefault}{\mddefault}{\updefault}{\color[rgb]{0,0,0}$r'$}%
}}}
\put(2160,-3091){\makebox(0,0)[lb]{\smash{\SetFigFont{5}{6.0}{\rmdefault}{\mddefault}{\updefault}{\color[rgb]{0,0,0}$\beta < 0$}%
}}}
\put(5634,367){\makebox(0,0)[lb]{\smash{\SetFigFont{5}{6.0}{\rmdefault}{\mddefault}{\updefault}{\color[rgb]{0,0,0}$R$}%
}}}
\put(5655,-1748){\makebox(0,0)[lb]{\smash{\SetFigFont{5}{6.0}{\rmdefault}{\mddefault}{\updefault}{\color[rgb]{0,0,0}$r'$}%
}}}
\put(7410,-3091){\makebox(0,0)[lb]{\smash{\SetFigFont{5}{6.0}{\rmdefault}{\mddefault}{\updefault}{\color[rgb]{0,0,0}$\beta > 0$}%
}}}
\put(9969,-2334){\makebox(0,0)[lb]{\smash{\SetFigFont{5}{6.0}{\rmdefault}{\mddefault}{\updefault}{\color[rgb]{0,0,0}$K$}%
}}}
\put(8490,-2265){\makebox(0,0)[lb]{\smash{\SetFigFont{5}{6.0}{\rmdefault}{\mddefault}{\updefault}{\color[rgb]{0,0,0}$k$}%
}}}
\put(7868,-2280){\makebox(0,0)[lb]{\smash{\SetFigFont{5}{6.0}{\rmdefault}{\mddefault}{\updefault}{\color[rgb]{0,0,0}$k'$}%
}}}
\put(5692,-1148){\makebox(0,0)[lb]{\smash{\SetFigFont{5}{6.0}{\rmdefault}{\mddefault}{\updefault}{\color[rgb]{0,0,0}$r$}%
}}}
\end{picture}

   \caption{Geography in the case (b).}
   \label{fig:case-b}
\end{center}
\end{figure}

Consider the correspondence $\cA_{r,r'}^{\chi(v),k'}$,
\begin{equation}
   \correspondence{BN_{k}^o(A_l(v))}{BN_{\chi(v)}^o(v)}{BN_{k'}^0(v')}{\pi_1}{\pi_2}
   \label{diagram:case-b}
\end{equation}

as illustrated by the Figure ~\ref{fig:case-b}.
  We can assume
$$
   \rho(\br,g,d) \ge 0,
$$
i.e., $(k,r)$ and $(k',r') \in D_1$, since otherwise otherwise
$M(r,a,\beta)$ and all the varieties in the correspondence
~\ref{diagram:case-b} are empty.

  Since $r-r' = k-k'$ and $k = 2r+\beta$, we have $r=k'-r' - \beta$
and $l = r-r' = k'- 2r' - \beta$.

   Now, if we are given some $v'=(r',a,\beta)$ and $k' \ge \chi(v') =
2r'+\beta$, we define $r$ by the formula given above, we can include
$BN_{k'}(v')$ in the correspondence of the form ~\ref{diagram:case-b}.

\begin{thing}
Proof of the main theorem:

We start by proving part 1. Let $(k,r) \in D_1$, and assume that
$(V_d^{\br})''(|a|^s)$ is not empty. Let $(C,B) \in
(V_d^{\br})''(|a|^s)$. By Lemma ~\ref{conclusion-of-case-a}, the
scheme $BN_{\chi(v)}^0(M(v))$ is not empty and contains a globally
generated vector bundle. Now, by Lemmas
~\ref{generic-evaluation-maps-in-globally-generated-case-are-monomorphic},
~\ref{generic-factors-in-globally-generated-case-are-acceptable} and
~\ref{stability-of-factors} the variety $BN_k^o(A_l(v))$ is not empty,
and thus $BN_{k'}(v')$ is not empty.

Now the variety $BN_k^o(v)$ is a nonempty open subscheme in $M(v)$,
and thus 
$$
    \dim BN_{\chi(v)}^0(v) = 
    \dim M(v) = 2 \rho(\br,g,d)
$$

Let $BN_{\chi(v)}^{o,gg}(v)$ be an open subscheme in
$BN_{\chi(v)}^{o,gg}(v)$ parametrizing globally generated sheaves; we
have just proved that it is not empty.  The lemmas
~\ref{generic-evaluation-maps-in-globally-generated-case-are-monomorphic},
~\ref{generic-factors-in-globally-generated-case-are-acceptable}, and
~\ref{stability-of-factors} imply that the fiber of $\pi_1$ over a
point of $BN_{\chi(v)}^{o,gg}(v)$ is isomorphic to a nonempty open
subscheme in the Grassmanian variety $Gr(l,2r+\beta)$ of dimension
$l(2r+\beta-l)$.  (Note that since $l \le 2r+\beta$, this number is
not negative). In particular, $\pi_1$ is epimorphic over
$BN_{\chi(v)}^{o,gg}(v)$.

% Fibers of $\pi_2$:
% For the extension to be stable 
% we need injectivity of $\alpha_e$ and 
% either $(c_1(F),H) = 1$, or $\Pic S = Z$ (the second one is in the 
% condition of the Main Theorem)
%
% For the extension to be acceptable we need
%     nothing, if $r' >= 2$;
%     Caley-Bacharash,  if  $r' = 1'$;
%     nothing, if $r' = 0$, and $\dim W = 1$;
%     $A$ is globally generated,  $r' = 0$, $\dim W >= 2$.
%  

Let us study the fibers of $\pi_2$. The lemmas
~\ref{generic-extensions-are-acceptable} and
~\ref{stability-of-extensions} imply that

%if either one of the following
%
%\begin{itemize}
%    \item[] $r' \ge 2$;
%    \item[] $r' = 1$, and [$E'$ is Caley-Bacharash];
%    \item[] $r' = 0$, and $l = 1$;
%    \item[] $r' = 0$, and $l \ge 2$, $E' = (i_C)_*(B)$, and
%                              $A=K_C B^{-1}$ is globally generated
%\end{itemize}
%%
%holds, then $\pi_2^{-1}(E')$ is not empty. It follows that

%
\begin{enumerate}
    \item If $r' \ge 2$, then $\pi_2$ is epimorphic;
    \item If $r' = 1$, then $\pi_2$ is epimorphic over a
               subscheme of simple Caley-Bacharash 0-cycles;
    \item If $r' = 0$ and $l = 1$, then $\pi_2$ is epimorphic, and
    \item If $r' = 0$ and $l \ge 2$, then  $\pi_2$ is 
          epimorphic over $  (V_{2g-2-d}^{\br + \beta})'(|a|^s)$.
\end{enumerate}

For the proof of part one of the main theorem, we need to consider
only the case $r' \ge 2$, when $\pi_2$ is epimorphic.

The fiber of $\pi_2$ over a point $E' \in BN_{k'}(M(v'))$ is
isomorphic to an open subscheme in the Grassmanian variety
$\Gr(l,\Ext^1_S(E',\os)) \simeq \Gr(l,H^1(S,E')\dual) \simeq \Gr(l,
h^1(S,E'))$. It is easy to see geometrically that $h^1(S,E') = l$, 
and it follows that $\pi_2$ is an isomorphism. It follows
that the correspondence $\cA$ in the case (b) gives
a morphism
$$
   \pi_2^{-1}\pi_1: \; \; 
   BN_k^o(r',a,\beta) \to BN_{2r+\beta}^0(r,a,\beta),
$$
where $r=k'-r' - \beta$, as denoted with an arrow on the Figure
~\ref{fig:case-b}. 

 Now counting of dimensions gives
\begin{multline*}
  \dim BN_{k'}(M(v')) = 
  \dim M(v) + \dim \text{(fiber of $\pi_1$) } = \\ =
  2(g-r(r-\alpha)) + l(2r+\beta-l) = \\ =
  2(g-r'(r'-\alpha)) - k'(k'- (2r'+\beta) ) = \dim M(v') -
  k'(k'-\chi(v')) = \\ = \vdim BN_{k'}(v').
\end{multline*}

The proof of the part 0 of the main theorem is parallel, except that
our study in the case (a) in this case does not guaranty that
$BN_k(A_l^o(v))$ is not empty.

To prove parts 2 and 3 of the main theorem, assume that $(k',r') \in
D_2 \cup D_3$. Assume that $BN_k'(v')$ is not empty. Then, by Lemmas
~\ref{generic-extensions-are-acceptable} and
~\ref{stability-of-extensions}, the varieties $BN_k^0(A_l(v))$ and
thus $BN_k^o(v)$ are not empty, which contradicts to the Mukai
theorem, since whenever $M(r,a,\beta)$ is not empty, the dimension
$\dim M(r,a,\beta) = 2 \rho (\br,g,d)$ is not negative.

%-------------------------------------------------------------------
%
\subsubsection{Case (c)}
\label{section:case-c}

Let us fix $v=(r,a,\beta) \in K'(S)$ and $k \in \ZZ$ such that $r \ge
2$ and $k \ge \max(\chi(v),r) = \max(2r-\beta,r)$.  Consider the
correspondence $\cA_{r,0}^{k,k-r}$,
$$
   \correspondence{BN_{k}^o(A_{r}(v))}{BN_{k}^o(v)}{%
         V_{2g-2-d}^{k-\br}(|a|) }{\pi_1^k}{\pi_2^k}
$$

%The fibers of $\pi_1^k$ are open subschemes in the Grassmanian
%varieties $\Gr(r,k)$, and the fibers of $\pi_2^k$ are open subschemes
%in the $\Gr(r-1,k-r-\beta)$.

Let $(k_0,r_0)$ be the intersection point of the lines $l(\beta)$ and
$l'(k,r)$, as on the Figure ~\ref{fig:case-c}.
\begin{figure}[!h]
\begin{center}
\setlength{\unitlength}{3947sp}%
\begingroup\makeatletter\ifx\SetFigFont\undefined%
\gdef\SetFigFont#1#2#3#4#5{%
  \reset@font\fontsize{#1}{#2pt}%
  \fontfamily{#3}\fontseries{#4}\fontshape{#5}%
  \selectfont}%
\fi\endgroup%
\begin{picture}(1756,1402)(732,-644)
{\color[rgb]{0,0,0}\thinlines
\put(2101,389){\circle{12}}
}%
{\color[rgb]{0,0,0}\put(1726, 14){\circle{12}}
}%
{\color[rgb]{0,0,0}\put(1201,-511){\circle{12}}
}%
{\color[rgb]{0,0,0}\put(901,-211){\line( 2, 1){1500}}
}%
{\color[rgb]{0,0,0}\multiput(2101,389)(0.00000,-85.71429){11}{\line( 0,-1){ 42.857}}
}%
{\color[rgb]{0,0,0}\put(901,689){\line( 0,-1){1200}}
\put(901,-511){\line( 1, 0){1575}}
}%
{\color[rgb]{0,0,0}\multiput(901,389)(88.88889,0.00000){14}{\line( 1, 0){ 44.444}}
}%
{\color[rgb]{0,0,0}\put(1726, 14){\line(-1,-1){525}}
}%
{\color[rgb]{0,0,0}\multiput(901, 14)(126.92308,0.00000){7}{\line( 1, 0){ 63.462}}
}%
{\color[rgb]{0,0,0}\multiput(1726, 14)(6.35593,6.35593){60}{\makebox(1.6667,11.6667){\SetFigFont{5}{6}{\rmdefault}{\mddefault}{\updefault}.}}
}%
{\color[rgb]{0,0,0}\multiput(1726, 14)(0.00000,-95.45455){6}{\line( 0,-1){ 47.727}}
}%
\put(2479,-604){\makebox(0,0)[lb]{\smash{\SetFigFont{5}{6.0}{\familydefault}{\mddefault}{\updefault}{\color[rgb]{0,0,0}$K$}%
}}}
\put(789,650){\makebox(0,0)[lb]{\smash{\SetFigFont{5}{6.0}{\familydefault}{\mddefault}{\updefault}{\color[rgb]{0,0,0}$R$}%
}}}
\put(2009,-588){\makebox(0,0)[lb]{\smash{\SetFigFont{5}{6.0}{\familydefault}{\mddefault}{\updefault}{\color[rgb]{0,0,0}$k_0$}%
}}}
\put(1654,-586){\makebox(0,0)[lb]{\smash{\SetFigFont{5}{6.0}{\familydefault}{\mddefault}{\updefault}{\color[rgb]{0,0,0}$k$}%
}}}
\put(732,364){\makebox(0,0)[lb]{\smash{\SetFigFont{5}{6.0}{\familydefault}{\mddefault}{\updefault}{\color[rgb]{0,0,0}$r_0$}%
}}}
\put(799,-14){\makebox(0,0)[lb]{\smash{\SetFigFont{5}{6.0}{\familydefault}{\mddefault}{\updefault}{\color[rgb]{0,0,0}$r$}%
}}}
\put(785,-221){\makebox(0,0)[lb]{\smash{\SetFigFont{5}{6.0}{\familydefault}{\mddefault}{\updefault}{\color[rgb]{0,0,0}$\frac{\alpha}{2}$}%
}}}
\put(2176,539){\makebox(0,0)[lb]{\smash{\SetFigFont{5}{6.0}{\familydefault}{\mddefault}{\updefault}{\color[rgb]{0,0,0}$l$}%
}}}
\put(1471,-124){\makebox(0,0)[lb]{\smash{\SetFigFont{5}{6.0}{\familydefault}{\mddefault}{\updefault}{\color[rgb]{0,0,0}$l'$}%
}}}
\end{picture}

   \caption{Geography in the case (c).}
   \label{fig:case-c}
\end{center}
\end{figure}
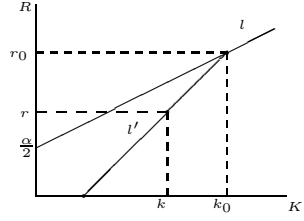
We have $r_0 = k - r - \beta$. The results of
~\ref{section:case-a} and ~\ref{section:case-b}
%the cases (a) and (b) 
ensure that if
\begin{enumerate}
   \item $(k,r) \in D_1$,
   \item The variety $(V_d^{\br_0}(|a|^s))''$ is not empty, and
   \item $\Pic S = \ZZ a$ or $(a,H) = 1$,
\end{enumerate}
then all the terms in the correspondence $\cA_{r,0}^{k,k-r}$ are not
empty and the morphisms $\pi_1^k$ and $\pi_2^k$ are dominant.

Let us fix $v$ and vary $k \in \ZZ$, $k \ge \chi(v)$.
As one can see from the Figure ~\ref{fig:r-0},
\begin{figure}[!h]
\begin{center}
\setlength{\unitlength}{3522sp}%
\begingroup\makeatletter\ifx\SetFigFont\undefined%
\gdef\SetFigFont#1#2#3#4#5{%
  \reset@font\fontsize{#1}{#2pt}%
  \fontfamily{#3}\fontseries{#4}\fontshape{#5}%
  \selectfont}%
\fi\endgroup%
\begin{picture}(5715,1873)(328,-1027)
\thinlines
{\color[rgb]{0,0,0}\put(451,-106){\line( 2, 1){1170}}
}%
{\color[rgb]{0,0,0}\put(946, 74){\vector(-1,-1){495}}
}%
{\color[rgb]{0,0,0}\put(2521,-106){\line( 2, 1){1350}}
}%
{\color[rgb]{0,0,0}\put(2521,-421){\vector(-1,-1){  0}}
\put(2521,-421){\vector( 1, 1){630}}
}%
{\color[rgb]{0,0,0}\put(4771,-421){\vector( 1, 1){900}}
}%
{\color[rgb]{0,0,0}\put(451,749){\line( 0,-1){1170}}
\put(451,-421){\line( 1, 0){1440}}
}%
{\color[rgb]{0,0,0}\multiput(1081,209)(-140.00000,0.00000){5}{\line(-1, 0){ 70.000}}
}%
{\color[rgb]{0,0,0}\multiput(1081,209)(0.00000,-140.00000){5}{\line( 0,-1){ 70.000}}
}%
{\color[rgb]{0,0,0}\multiput(451, 74)(12.00000,0.00000){121}{\makebox(1.8676,13.0735){\SetFigFont{5}{6}{\rmdefault}{\mddefault}{\updefault}.}}
}%
{\color[rgb]{0,0,0}\multiput(946, 74)(90.00000,90.00000){2}{\line( 1, 1){ 45.000}}
}%
{\color[rgb]{0,0,0}\multiput(946, 74)(0.00000,-141.42857){4}{\line( 0,-1){ 70.714}}
}%
{\color[rgb]{0,0,0}\put(1756, 74){\line(-1,-1){495}}
}%
{\color[rgb]{0,0,0}\put(2521,-421){\line( 1, 0){1440}}
}%
{\color[rgb]{0,0,0}\put(2521,749){\line( 0,-1){1170}}
}%
{\color[rgb]{0,0,0}\multiput(2521,209)(140.00000,0.00000){5}{\line( 1, 0){ 70.000}}
}%
{\color[rgb]{0,0,0}\multiput(3151,209)(9.00000,0.00000){81}{\makebox(1.8676,13.0735){\SetFigFont{5}{6}{\rmdefault}{\mddefault}{\updefault}.}}
}%
{\color[rgb]{0,0,0}\put(3691,209){\line(-1,-1){630}}
}%
{\color[rgb]{0,0,0}\multiput(3151,209)(0.00000,-140.00000){5}{\line( 0,-1){ 70.000}}
}%
{\color[rgb]{0,0,0}\put(4501,749){\line( 0,-1){1170}}
}%
{\color[rgb]{0,0,0}\put(4501,-421){\line( 1, 0){1530}}
}%
{\color[rgb]{0,0,0}\put(4501,-106){\line( 2, 1){1530}}
}%
{\color[rgb]{0,0,0}\multiput(4501,-421)(96.92308,96.92308){7}{\line( 1, 1){ 48.462}}
}%
{\color[rgb]{0,0,0}\multiput(5131,209)(0.00000,-140.00000){5}{\line( 0,-1){ 70.000}}
}%
{\color[rgb]{0,0,0}\multiput(5131,209)(-140.00000,0.00000){5}{\line(-1, 0){ 70.000}}
}%
{\color[rgb]{0,0,0}\multiput(5671,479)(-137.64706,0.00000){9}{\line(-1, 0){ 68.824}}
}%
{\color[rgb]{0,0,0}\multiput(5671,479)(9.00000,0.00000){41}{\makebox(1.8676,13.0735){\SetFigFont{5}{6}{\rmdefault}{\mddefault}{\updefault}.}}
}%
{\color[rgb]{0,0,0}\put(5941,479){\line(-1,-1){900}}
}%
\put(351,738){\makebox(0,0)[lb]{\smash{\SetFigFont{5}{6.0}{\familydefault}{\mddefault}{\updefault}{\color[rgb]{0,0,0}$R$}%
}}}
\put(3077,-928){\makebox(0,0)[lb]{\smash{\SetFigFont{5}{6.0}{\familydefault}{\mddefault}{\updefault}{\color[rgb]{0,0,0}$r = \alpha$}%
}}}
\put(5035,-489){\makebox(0,0)[lb]{\smash{\SetFigFont{5}{6.0}{\familydefault}{\mddefault}{\updefault}{\color[rgb]{0,0,0}$\alpha$}%
}}}
\put(5434,-487){\makebox(0,0)[lb]{\smash{\SetFigFont{5}{6.0}{\familydefault}{\mddefault}{\updefault}{\color[rgb]{0,0,0}$2r-\alpha$}%
}}}
\put(4605,-481){\makebox(0,0)[lb]{\smash{\SetFigFont{5}{6.0}{\familydefault}{\mddefault}{\updefault}{\color[rgb]{0,0,0}$r-\alpha$}%
}}}
\put(4864,-978){\makebox(0,0)[lb]{\smash{\SetFigFont{5}{6.0}{\familydefault}{\mddefault}{\updefault}{\color[rgb]{0,0,0}$r > \alpha$}%
}}}
\put(2422,-113){\makebox(0,0)[lb]{\smash{\SetFigFont{5}{6.0}{\familydefault}{\mddefault}{\updefault}{\color[rgb]{0,0,0}$\frac{\alpha}{2}$}%
}}}
\put(3587,478){\makebox(0,0)[lb]{\smash{\SetFigFont{5}{6.0}{\familydefault}{\mddefault}{\updefault}{\color[rgb]{0,0,0}$l$}%
}}}
\put(872,-978){\makebox(0,0)[lb]{\smash{\SetFigFont{5}{6.0}{\familydefault}{\mddefault}{\updefault}{\color[rgb]{0,0,0}$r < \alpha$}%
}}}
\put(2415,700){\makebox(0,0)[lb]{\smash{\SetFigFont{5}{6.0}{\familydefault}{\mddefault}{\updefault}{\color[rgb]{0,0,0}$R$}%
}}}
\put(4376,701){\makebox(0,0)[lb]{\smash{\SetFigFont{5}{6.0}{\familydefault}{\mddefault}{\updefault}{\color[rgb]{0,0,0}$R$}%
}}}
\put(897,-486){\makebox(0,0)[lb]{\smash{\SetFigFont{5}{6.0}{\familydefault}{\mddefault}{\updefault}{\color[rgb]{0,0,0}$r$}%
}}}
\put(983,-494){\makebox(0,0)[lb]{\smash{\SetFigFont{5}{6.0}{\familydefault}{\mddefault}{\updefault}{\color[rgb]{0,0,0}$\alpha$}%
}}}
\put(1912,-518){\makebox(0,0)[lb]{\smash{\SetFigFont{5}{6.0}{\familydefault}{\mddefault}{\updefault}{\color[rgb]{0,0,0}$K$}%
}}}
\put(1354,413){\makebox(0,0)[lb]{\smash{\SetFigFont{5}{6.0}{\familydefault}{\mddefault}{\updefault}{\color[rgb]{0,0,0}$l$}%
}}}
\put(3020,-499){\makebox(0,0)[lb]{\smash{\SetFigFont{5}{6.0}{\familydefault}{\mddefault}{\updefault}{\color[rgb]{0,0,0}$\alpha$}%
}}}
\put(3900,-496){\makebox(0,0)[lb]{\smash{\SetFigFont{5}{6.0}{\familydefault}{\mddefault}{\updefault}{\color[rgb]{0,0,0}$K$}%
}}}
\put(352,182){\makebox(0,0)[lb]{\smash{\SetFigFont{5}{6.0}{\familydefault}{\mddefault}{\updefault}{\color[rgb]{0,0,0}$\alpha$}%
}}}
\put(328,-128){\makebox(0,0)[lb]{\smash{\SetFigFont{5}{6.0}{\familydefault}{\mddefault}{\updefault}{\color[rgb]{0,0,0}$\frac{\alpha}{2}$}%
}}}
\put(355, 46){\makebox(0,0)[lb]{\smash{\SetFigFont{5}{6.0}{\familydefault}{\mddefault}{\updefault}{\color[rgb]{0,0,0}$r$}%
}}}
\put(2431,192){\makebox(0,0)[lb]{\smash{\SetFigFont{5}{6.0}{\familydefault}{\mddefault}{\updefault}{\color[rgb]{0,0,0}$\alpha$}%
}}}
\put(4401,451){\makebox(0,0)[lb]{\smash{\SetFigFont{5}{6.0}{\familydefault}{\mddefault}{\updefault}{\color[rgb]{0,0,0}$r$}%
}}}
\put(4376,-127){\makebox(0,0)[lb]{\smash{\SetFigFont{5}{6.0}{\familydefault}{\mddefault}{\updefault}{\color[rgb]{0,0,0}$\frac{\alpha}{2}$}%
}}}
\put(5998,-514){\makebox(0,0)[lb]{\smash{\SetFigFont{5}{6.0}{\familydefault}{\mddefault}{\updefault}{\color[rgb]{0,0,0}$K$}%
}}}
\put(5738,583){\makebox(0,0)[lb]{\smash{\SetFigFont{5}{6.0}{\familydefault}{\mddefault}{\updefault}{\color[rgb]{0,0,0}$l$}%
}}}
\end{picture}

   \caption{Geography in the case (c), three variants}
   \label{fig:r-0}
\end{center}
\end{figure}
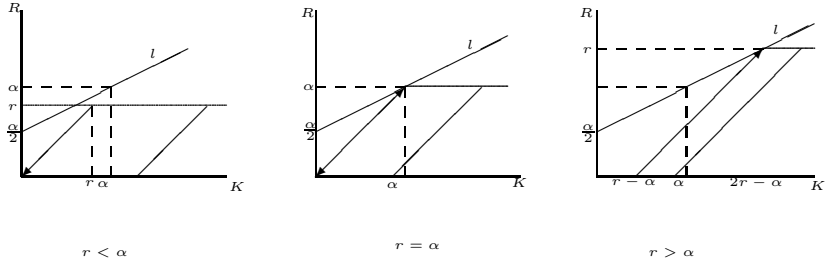
there are three combinatorial variants, when $r < \alpha$, $r =
\alpha$ and $r > \alpha$.

% r < \alpha; d > g+r-1
%
  If $r < \alpha$, then the tom term in the correspondence
$\cA_{r,0}^{k,k-r}$ is empty for $k < r$, and all the terms are
nonempty if $k \ge r$. Additionally, if $k=r$, then $\pi_1^r$ is an
isomorphism over a subscheme of globally generated vector bundles in
$BN_{r}(v)$, and thus $\pi_1^r$ is a birational isomorphism and the
correspondence $\cA_{r,0}^{r,0}$ gives a rational map
$BN_{r}(r,a,\beta) \to \pic^d(|a|^s)$, marked with an arrow on the
Figure ~\ref{fig:r-0}.

% r = \alpha; d = g+r-1
%
  If $r = \alpha$, the correspondence $\cA_{r,0}^{r,0}$ gives a
birational isomorphism
$$
   M(r,a,-r) \stackrel{bir}{\simeq} \pic^{2g-2-d}(|a|)
$$
where $d = g+ r = \frac{a^2+2}{2} + r$.

% r > \alpha; d < g+r-1
%
  If $r > \alpha$, the map $\pi_2^k$ is a birational isomorphism for
$k = \chi(v)$, and the correspondence $\cA_{r,0}^{2r-\alpha,r-\alpha}$
gives a rational morphism $V_{2g-2-d}^{\br-\alpha}(|a|) \to M(v)$,
which is denoted by an arrow on the Figure ~\ref{fig:r-0}.

  It follows that 
\begin{enumerate}
   \item If $r < \alpha$ then the moduli space $M(v)$ contains a
Brill-Noether subscheme $BN_r(v)$ which is birationally isomorphic to
a Grassmanian fibration over the relative Picard scheme
$\pic^{2g-2-d}(|a|)$;
   \item If $r = \alpha$ then the moduli space $M(v)$ is birationally
isomorphic to the relative Picard variety $\pic^{2g-2-d}(|a|)$, and
   \item If $r > \alpha$, then the relative Picard variety
$\pic^{2g-2-d}(|a|)$ contains a relative Brill-Noether subscheme
$W_{2g-2-d}^{\br-\alpha}(|a|)$ which is birationally isomorphic to a
Grassmanian fibration over the moduli space $M(v)$.
\end{enumerate}
\subsubsection{Case (d)}
\label{section:case-d}

Let us fix $v=(r,a,\beta)$ such that $r \ge 2$ and $k \ge
\max(\chi(v),r-1)$. The fibers of the map $\pi_1^k$ in the
correspondence $\cA_{r,1}^{k,k-(r-1)}$
$$
   \correspondence{BN_{k}^o(A_{r-1}(v))}{BN_{k}^o(v)}{
   \hilb^d_{(L,k-\br - \beta)^o}(S) }{\pi_1^k}{\pi_2^k}
$$
are isomorphic to open subschemes in the Grassmanian $\Gr(r-1,k)$, and
the fibers of $\pi_2^k$ are isomorphic to open subschemes in the
Grassmanian $\Gr(r-1,k-r+\alpha-1)$. (Note that $r-1 \le k-r+\alpha-1$,
and $r-1 = k-r+\alpha-1$ only if  $k = \chi(v)$).

Let $(k_0,r_0)$ be the intersection point of $l(\beta)$ and $l'(k,r)$.
We have $r_0 = k - r - \beta$. Lemma ~\ref{conclusion-of-case-a}
implies that if
\begin{enumerate}
   \item $(k,r) \in D_1$,
   \item The variety $(V_d^{\br_0}(|a|^s))''$ is not empty, and
   \item $\Pic S = \ZZ a$ or $(a,H) = 1$,
\end{enumerate}
then $BN_{\chi_0}(v_0)$ is not empty and contains a globally generated
vector bundle, and Lemmas
~\ref{generic-evaluation-maps-in-globally-generated-case-are-monomorphic}
and ~\ref{generic-factors-in-globally-generated-case-are-acceptable}
imply that $\hilb^d_{(L,k-\br - \beta)^o}(S)$ is not empty. One can
see that the morphisms $\pi_1^k$ and $\pi_2^k$ are dominant.

% thus we do not use the nonemptiness criteria for $sCB()$
% which we proved independently.

There are three combinatorial possibilities, depending on whether $r <
\alpha -1$, $r = \alpha -1$ ($d=r+g$) or $r > \alpha - 1$., as one can
see on the Figure ~\ref{fig:r-1.3-cases}.
%
% demonstrates the geometry of
%   he correspondence $\cA_{r,1}^{k,k-(r-1)}$ in these three cases.
%

\begin{figure}[!h]
\begin{center}
\setlength{\unitlength}{3522sp}%
\begingroup\makeatletter\ifx\SetFigFont\undefined%
\gdef\SetFigFont#1#2#3#4#5{%
  \reset@font\fontsize{#1}{#2pt}%
  \fontfamily{#3}\fontseries{#4}\fontshape{#5}%
  \selectfont}%
\fi\endgroup%
\begin{picture}(5580,1758)(95,-916)
\thinlines
{\color[rgb]{0,0,0}\put(361,749){\line( 0,-1){1080}}
\put(361,-331){\line( 1, 0){1260}}
}%
{\color[rgb]{0,0,0}\put(361, 74){\line( 1, 0){1260}}
}%
{\color[rgb]{0,0,0}\put(401,-216){\makebox(1.8676,13.0735){\SetFigFont{5}{6}{\rmdefault}{\mddefault}{\updefault}.}}
}%
{\color[rgb]{0,0,0}\put(676, 74){\vector(-1,-1){315}}
}%
{\color[rgb]{0,0,0}\put(2431,749){\line( 0,-1){1080}}
\put(2431,-331){\line( 1, 0){1170}}
}%
{\color[rgb]{0,0,0}\put(2431,-241){\line( 1, 0){1170}}
}%
{\color[rgb]{0,0,0}\put(2431,209){\line( 1, 0){1170}}
}%
{\color[rgb]{0,0,0}\put(2431,-241){\vector(-1,-1){  0}}
\put(2431,-241){\vector( 1, 1){450}}
}%
{\color[rgb]{0,0,0}\put(4411,749){\line( 0,-1){1080}}
\put(4411,-331){\line( 1, 0){1170}}
}%
{\color[rgb]{0,0,0}\put(4411,-16){\line( 2, 1){1170}}
}%
{\color[rgb]{0,0,0}\put(4411,-241){\line( 1, 0){1170}}
}%
{\color[rgb]{0,0,0}\multiput(2431,-16)(8.06897,4.03448){146}{\makebox(1.8676,13.0735){\SetFigFont{5}{6}{\rmdefault}{\mddefault}{\updefault}.}}
}%
{\color[rgb]{0,0,0}\multiput(676, 74)(90.00000,90.00000){2}{\line( 1, 1){ 45.000}}
}%
{\color[rgb]{0,0,0}\multiput(676, 74)(0.00000,-115.71429){4}{\line( 0,-1){ 57.857}}
}%
{\color[rgb]{0,0,0}\multiput(811,209)(-128.57143,0.00000){4}{\line(-1, 0){ 64.286}}
}%
{\color[rgb]{0,0,0}\multiput(811,209)(0.00000,-120.00000){5}{\line( 0,-1){ 60.000}}
}%
{\color[rgb]{0,0,0}\multiput(2881,209)(0.00000,-120.00000){5}{\line( 0,-1){ 60.000}}
}%
{\color[rgb]{0,0,0}\put(1441, 74){\line(-1,-1){315}}
}%
{\color[rgb]{0,0,0}\put(3061,-241){\line( 1, 1){450}}
}%
{\color[rgb]{0,0,0}\put(361,-241){\line( 1, 0){1260}}
}%
{\color[rgb]{0,0,0}\put(4591,-241){\vector( 1, 1){630}}
}%
{\color[rgb]{0,0,0}\multiput(4411,-241)(100.00000,100.00000){5}{\line( 1, 1){ 50.000}}
}%
{\color[rgb]{0,0,0}\put(4411,389){\line( 1, 0){1170}}
}%
{\color[rgb]{0,0,0}\put(5491,389){\line(-1,-1){630}}
}%
{\color[rgb]{0,0,0}\multiput(4861,209)(-128.57143,0.00000){4}{\line(-1, 0){ 64.286}}
}%
{\color[rgb]{0,0,0}\multiput(4861,209)(0.00000,-120.00000){5}{\line( 0,-1){ 60.000}}
}%
{\color[rgb]{0,0,0}\multiput(5221,389)(0.00000,-130.90909){6}{\line( 0,-1){ 65.455}}
}%
{\color[rgb]{0,0,0}\put(361,-16){\line( 2, 1){1170}}
}%
\put(274,725){\makebox(0,0)[lb]{\smash{\SetFigFont{5}{6.0}{\familydefault}{\mddefault}{\updefault}{\color[rgb]{0,0,0}$R$}%
}}}
\put(1663,-444){\makebox(0,0)[lb]{\smash{\SetFigFont{5}{6.0}{\familydefault}{\mddefault}{\updefault}{\color[rgb]{0,0,0}$K$}%
}}}
\put(2379,-244){\makebox(0,0)[lb]{\smash{\SetFigFont{5}{6.0}{\familydefault}{\mddefault}{\updefault}{\color[rgb]{0,0,0}$1$}%
}}}
\put(4236,734){\makebox(0,0)[lb]{\smash{\SetFigFont{5}{6.0}{\familydefault}{\mddefault}{\updefault}{\color[rgb]{0,0,0}$R$}%
}}}
\put(5675,-391){\makebox(0,0)[lb]{\smash{\SetFigFont{5}{6.0}{\familydefault}{\mddefault}{\updefault}{\color[rgb]{0,0,0}$K$}%
}}}
\put(4500,-411){\makebox(0,0)[lb]{\smash{\SetFigFont{5}{6.0}{\familydefault}{\mddefault}{\updefault}{\color[rgb]{0,0,0}$r+\beta+1$}%
}}}
\put(3405,478){\makebox(0,0)[lb]{\smash{\SetFigFont{5}{6.0}{\familydefault}{\mddefault}{\updefault}{\color[rgb]{0,0,0}$l$}%
}}}
\put(705,-864){\makebox(0,0)[lb]{\smash{\SetFigFont{7}{8.4}{\familydefault}{\mddefault}{\updefault}{\color[rgb]{0,0,0}$r < \alpha - 1$}%
}}}
\put(2732,-856){\makebox(0,0)[lb]{\smash{\SetFigFont{7}{8.4}{\familydefault}{\mddefault}{\updefault}{\color[rgb]{0,0,0}$r = \alpha - 1$}%
}}}
\put(2298,691){\makebox(0,0)[lb]{\smash{\SetFigFont{5}{6.0}{\familydefault}{\mddefault}{\updefault}{\color[rgb]{0,0,0}$R$}%
}}}
\put(3600,-460){\makebox(0,0)[lb]{\smash{\SetFigFont{5}{6.0}{\familydefault}{\mddefault}{\updefault}{\color[rgb]{0,0,0}$K$}%
}}}
\put(1300,505){\makebox(0,0)[lb]{\smash{\SetFigFont{5}{6.0}{\familydefault}{\mddefault}{\updefault}{\color[rgb]{0,0,0}$l$}%
}}}
\put(2329,-31){\makebox(0,0)[lb]{\smash{\SetFigFont{5}{6.0}{\familydefault}{\mddefault}{\updefault}{\color[rgb]{0,0,0}$\frac{\alpha}{2}$}%
}}}
\put(2120,191){\makebox(0,0)[lb]{\smash{\SetFigFont{5}{6.0}{\familydefault}{\mddefault}{\updefault}{\color[rgb]{0,0,0}$\alpha-1$}%
}}}
\put(2756,-393){\makebox(0,0)[lb]{\smash{\SetFigFont{5}{6.0}{\familydefault}{\mddefault}{\updefault}{\color[rgb]{0,0,0}$\alpha - 2$}%
}}}
\put(4664,-867){\makebox(0,0)[lb]{\smash{\SetFigFont{7}{8.4}{\familydefault}{\mddefault}{\updefault}{\color[rgb]{0,0,0}$r > \alpha - 1$}%
}}}
\put(4300,368){\makebox(0,0)[lb]{\smash{\SetFigFont{5}{6.0}{\rmdefault}{\mddefault}{\updefault}{\color[rgb]{0,0,0}$r$}%
}}}
\put(5076,-406){\makebox(0,0)[lb]{\smash{\SetFigFont{5}{6.0}{\familydefault}{\mddefault}{\updefault}{\color[rgb]{0,0,0}$2r+\beta$}%
}}}
\put(5396,543){\makebox(0,0)[lb]{\smash{\SetFigFont{5}{6.0}{\familydefault}{\mddefault}{\updefault}{\color[rgb]{0,0,0}$l$}%
}}}
\put(284, 53){\makebox(0,0)[lb]{\smash{\SetFigFont{5}{6.0}{\familydefault}{\mddefault}{\updefault}{\color[rgb]{0,0,0}$r$}%
}}}
\put(239,-27){\makebox(0,0)[lb]{\smash{\SetFigFont{5}{6.0}{\familydefault}{\mddefault}{\updefault}{\color[rgb]{0,0,0}$\frac{\alpha}{2}$}%
}}}
\put(286,-255){\makebox(0,0)[lb]{\smash{\SetFigFont{5}{6.0}{\rmdefault}{\mddefault}{\updefault}{\color[rgb]{0,0,0}$1$}%
}}}
\put( 95,191){\makebox(0,0)[lb]{\smash{\SetFigFont{5}{6.0}{\familydefault}{\mddefault}{\updefault}{\color[rgb]{0,0,0}$\alpha - 1$}%
}}}
\put(4271,-33){\makebox(0,0)[lb]{\smash{\SetFigFont{5}{6.0}{\familydefault}{\mddefault}{\updefault}{\color[rgb]{0,0,0}$\frac{\alpha}{2}$}%
}}}
\put(4191,192){\makebox(0,0)[lb]{\smash{\SetFigFont{5}{6.0}{\rmdefault}{\mddefault}{\updefault}{\color[rgb]{0,0,0}$\alpha-1$}%
}}}
\put(4310,-258){\makebox(0,0)[lb]{\smash{\SetFigFont{5}{6.0}{\rmdefault}{\mddefault}{\updefault}{\color[rgb]{0,0,0}$1$}%
}}}
\end{picture}

   \caption{Geography in the case (d), three variants}
   \label{fig:r-1.3-cases}
\end{center}
\end{figure}
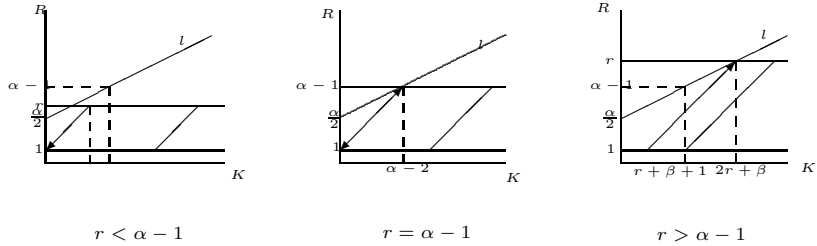
%

% r < \alpha -1; d > g+r
%
  In the first case $k \ge r-1$, the morphism $\pi_2^k$ has a positive
relative dimension and the morphism $\pi_1^k$ is a birational
isomorphism for $k = r-1$ and has a positive relative dimension for $k
> r-1$. It follows that for $k = r-1$ the correspondence $\cA_{r,1}^k$
gives a rational map
$$
   BN_{r-1}(r,a,\beta) \to \hilb^d(S),
$$
denoted with an arrow on the picture above.

% r = \alpha -1; d = g+r
%
If $r = \alpha -1$, then $k \ge r-1$ and for $k = r-1$ both $\pi_1^k$
and $\pi_2^k$ are birational isomorphisms, and thus we get the
birational isomorphism
$$
   M(r,a,-r-1) \stackrel{bir}{\simeq} \hilb^d (S)
$$
where $d = g+ r = \frac{a^2+2}{2} + r$.

[In the case $r = 2$ this statement was proved by Lothar Goettsche.]

% r > \alpha -1; d < g+r
%

If $r > \alpha -1$, then $k \ge \chi(v)$, the map $\pi_1^k$ has a
positive relative dimension, and $\pi_2^k$ is a birational isomorphism
for $k = 2r+ \beta$ and has a positive relative dimension for $k > 2r+
\beta$. It follows that for $k=2r+\beta$ ($\delta(\xi,L) = r-1$) the
correspondence gives a rational morphism $\hilb^d_{L,r-1}(S) \to M(v)$,
which is denoted by an arrow on the picture above.

Summarizing the results, we can see that for $r < \alpha - 1$ the
moduli space $M(v)$ contains a Brill-Noether subscheme which is
birational to a Grassmanian fibration over the Hilbert scheme
$\hilb^d(S)$, for $r = \alpha - 1$ the moduli space $M(v)$ is
birationally isomorphic to the Hilbert scheme $\hilb^d(S)$, and for $r
> \alpha - 1$ the Hilbert scheme $\hilb^d(S)$ contains a subscheme of
special 0-cycles which is birational to a Grassmanian fibration over
the moduli space $M(v)$.

\end{thing}
%-------------------------------------------------------------------

\subsection{Examples.}
\label{subsection:examples}

% ---------------------------------------------

\subsubsection{Example: $g=3, d=3$}
Let $(S,a)$ be a polarized K3 surface of genus $g = 3$, and let $d =
3$. Note that $M(0,a,\beta) = M(0,a,-1) \subset \Pic^{1}(|a|)$, since
that for $(C,B) \in M(0,a,\beta)$ we have $\deg B = 2g-2-d = 1$.

The formula ~\ref{eq:virtual-dimension} gives the following expected
dimensions of the Brill-Noether loci in various $M(r,a,-1)$:
$$
\begin{tabular}{c|ccccc}
     r &   &    &   &    & \\
     2 &   &    &   &  2 & \\
     1 &   &  6 &  4&  0 & \\
     0 &  6&  4 &  0&    & \\
  \hline
       &  0&  1 &  2&  3 & k 
\end{tabular}
$$

or, symbolically,
$$
\begin{tabular}{c|cccc}
     r &  &  &  &   \\ 
     2 &  &  &  & * \\ 
     1 &  & *& *& ! \\
     0 & *& *& !&   \\
  \hline
       & 0& 1& 2& 3 
\end{tabular}
$$
where we use ``*'' for integer points in $D_1$ and ``!'' for integer
points in $D_2$.

There are two integer points in $D_2$, corresponding to $V_1^1(|a|)$
and $\hilb^3_{(a,2)} S$. The variety $V_1^1(|a|)$ classifies linear
systems of type $g_1^1$ on smooth curves in $|a|$ and is evidently
empty for any $(S,a)$. The variety $\hilb^3_{(a,2)} S$ classifies
0-subschemes on $S$ of degree three with speciality index
$\delta(\xi,a) = 2$. Note that the linear system $|a|$ realizes $S$ as
a quartic in $\PP^3$ with at most double points. It follows that
$\hilb^3_{(a,2)} S$ is empty when $\phi_{|a|}$ is an embedding, and
parametrizes simple 0-cycles of the form $p+q+r$, where $p$, $q$ and
$r$ are three distinct points on a rational curve contracted by
$\phi_{|a|}$, and their degenerations. In particular, not that if
$\Pic S = \ZZ$, then $\hilb^3_{(a,2)} S$ is empty.

Now if $\Pic S = \ZZ$, then the main theorem guarantees that all the
varieties parametrized by points in $D_3$ are empty. 

\begin{figure}[!h]
\begin{center}
   \includegraphics[width=10cm,angle=-90]{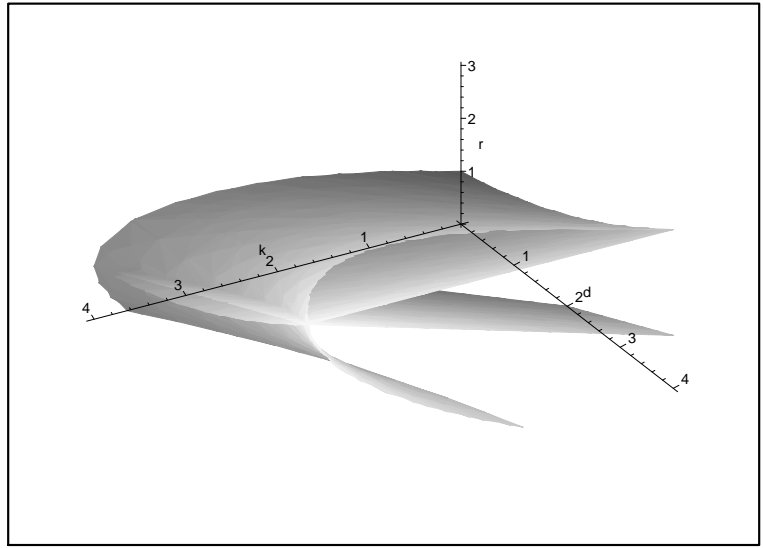}%
% n5   
   \caption{Geography for genus 3}
   \label{fig:geography.genus-3}
\end{center}
\end{figure}

One can consider the set of integer triples $(r,d,k)$ which satisfy
$$\vdim BN_k(r,a,\beta) \ge 0$$
%
% For $g = 3$ and $0 \le d \le 2g-2$ 
%
For $d \le 2g-2$ this domain is shown on the Figure %
~\ref{fig:geography.genus-3}.

%

% ------------------------------------------------------------------

\subsubsection{Example: $g=4, d=5$}
Let $(S,a)$ be a polarized K3 surface of genus $g = 4$, and let $d =
5$. Note that $M(0,a,\beta) = M(0,a,-2) \subset \Pic^{1}(|a|)$.  The
formula ~\ref{eq:virtual-dimension} gives the following table of the
Brill-Noether loci with nonnegative expected dimension,

\begin{comment}
%
$$
\begin{tabular}{c|cccccr}
   r &    &    &    &    &    &  \\
   3 &    &    &    &    &  2 &  \\
   2 &    &    &  8 &  5 &  0 &  \\
   1 & 10 &  9 &  6 &  1 &    &  \\
   0 &  8 &  5 &  0 &    &    &  \\
  \hline
     &  0 &  1 &  2 &  3 &  4 & k         
\end{tabular}
$$
or
\end{comment}

$$
\begin{tabular}{c|cccccc}
   r &   &   &   &   &   &  \\
   3 &   &   &   &   & * &  \\
   2 &   &   & * & * & ! &  \\
   1 & * & * & * & ! &   &  \\
   0 & * & * & ! &   &   &  \\ 
  \hline
     & 0 & 1 & 2 & 3 & 4 & k 
\end{tabular}
$$
where the notations are the same as in the previous example.

There are three Brill-Noether schemes in the region $D_2$.  For $r=0$
it is $V_1^1(|a|)$ which is empty. For $r=1$ it is $\hilb^5_{(a,3)}
S$.  A simple 0-cycle $\xi \in \hilb^5_{(a,3)}S$ is an intersection of
$S$ with a 3-secant line $l$ in the projective embedding given by
$|a|$. There is a net of hyperplanes in $\PP^4$ through $l$, and each
hyperplane cuts out a genus 4 curve on $S$ containing $\xi$. On such a
curve $\xi$ gives a linear system of type $g_5^3$. Since a genus 5
curve does not have any $g_5^3$'s, this variety is empty. (Essentially
this proof, which is a reduction from the case $r=1$ to $r=0$, can be
described as a correspondence between the Brill-Noether schemes
$V_{2g-2-d}^r(|a|^s)$ and $\Hilb^{d}_{(r-\beta,a)}(S)$).

  If $r=2$ and $\Pic S = \ZZ a$, then the emptiness of $BN_4(2,a,-2)$
follows from the discussion in ~\ref{section:case-b}.

%--------------------------- ------------------------------------
%
\subsubsection{Geometry of $W_4^1(C)$ for a trigonal curve of genus five} 
(Compare ~\cite{ACGH}).
\label{example:two-components-of-W_4^1(C)}

We will need the results of this example in the examples
~\ref{example:g=5,d=4} and ~\ref{example:g=5,d=5}. For a slightly
different approach, see ~\cite{ACGH}.

Let $C$ be a trigonal curve of genus five. It is well known that $C$
has unique trigonal structure $g_3^1$. Let $t \in |g_3^1|$ be any
divisor consisting of three distinct points. It follows that the
linear span $<t>$ is a trisecant line for $C$ in the canonical
embedding of $C$. We denote this line as $l$. 

Since $K_C - g_3^1 = g_5^2$, $C$ has a (unique) structure of a plane
quintic. Note that this linear system can be realized as a linear
system of hyperplanes through $l$ in the canonical embedding of
$C$. Let $C'$ be the plane model of $C$.  Since genus of a smooth
plane quintic is equal to $6$, it follows that the of $C'$ has one
double point $q'$, which can be either node or cusp.

Consider first the case when $q'$ is a node, and let $q_1$ and $q_2$
be two preimages of $q'$ on $C$. Since $K_C - g_3^1$ does not separate
$q_1$ and $q_2$, we have $<l + q_1> = <l + q_2>$ for the linear spans
of divisors. Let $\pi = <l+q_i>$; it is a two-dimensional plane in
$\PP^4$. 

Note that the divisor $t+q_1+q_2$ is of type $g_5^2$. Let $H$ be a
hyperplane in $\PP^4$ through $\pi$. It intersects $C$ in a degree 3
divisor $t'$ which forms $g_3^1$.  Since $|K_C| = |t + t' + q_1 +
q_2|$, It follows that $q_1 + q_2 \in |K_C - 2t|$.

Let $p$ be any point on $C$. We define $L_p = g_5^2 - p$, $M_p = g_3^1
+ p$. Note that the Serre duality interchanges $L_p$ and $M_p$: 
$$
   K_C L_p^{-1} \simeq M_p
$$

We give the following geometric interpretation of $|L_p|$. Note that
for any $p \in C$ we have $p \notin l$, since there are no linear
systems of type $g_4^2$ on $C$. (Note that the definition of the linear
span $<t+p>$ includes an infinitesimal computation if $p \in t$.)  It
follows that $pi_p=<l+p>$ is a two-plane in $\PP^4$. Projection with
center $\pi_p$ gives a linear system of type $g_4^1$ on $C$ which is
equal to $|K_C - t - p| = |g_5^2 - p| = |L_p|$.  (In terms of the
plane model $C'$, $|L_p|$ is a projection from the point $p$ on $C'$).

If $p \ne q_i$, then $|L_p|$ forms a $g_4^1$ without fixed points.
Note that $q_2$ is a fixed point of $|L_{q_1}|$ and $q_1$ is a fixed
point of $|L_{q_2}|$.

   Note that 
\begin{gather*}
   L_{q_1} = |K_C - t - q_1| = |t' + q_2| = M_{q_2},  \\
   L_{q_2} = M_{q_1}
\end{gather*}
   In particular, $K_C L_{q1}^{-1} = M_{q_1} = L_{q_2}$.

{\bf Lemma.} Any linear system of type $g_4^1$ on $C$ is of type $L_p$
or $M_p$ for some $p \in C$.

{\bf Proof.}  Let $L$ be a linear system of type $g_4^1$, and let $D
\in |L|$ be a divisor consisting of 4 distinct points on $C$: $D = d_1
+ d_2 + d_3 + d_4$. Let $\pi = <D>$; it is a 2-plane in $\PP^4$.

   Let $t \in |g_3^1|$ be a divisor consisting of three distinct
points, $t = t_1 + t_2 + t_3$. The linear span $l = <t>$ is a
trisecant line to $C$ in the canonical embedding.

   Case 1: $l \in \pi$. Consider the set-theoretic union $E = D \cup
t$. This is a set of 7, 6, 5 or 4 distinct points on $C$, and,
correspondingly, $|E|$ is a linear system of type $g_7^4$, $g_6^3$,
$g_5^2$ or $g_4^1$. Since $C$ does not have $g_7^4$ or $g_6^3$, the
order of $E$ is either 4 or 5.

  If $|E| = g_4^1$, then we have $t \subset D$, i.e., $|D| = |t + p| =
M_p$, and the lemma is proved.  Let $|E| = g_5^2$. After re-numeration
we have $D = t_1+t_2+d_3+d_4$. Note that the linear system
$g_5^2=K_C-t$ does not separate $d_3$ and $d_4$, since $<t,d_3> =
<t,d_4> = \pi$. It follows that, after renumeration, $d_3 = q_1$ and
$d_4 = q_2$. In particular,
$$
    |D| = |t_1 + t_2 + q_1 + q_2 | %
        = |t - t_3 + q_1 + q_2| = |g_5^2 - t_3| = L_{t_3},
$$
and the statement of the Lemma follows. 

   Case 2: $l \cap \pi = a$, where $a$ is a point in $\PP^4$.  Assume
first that $a = t_i$ for some $i$, for example, $a = t_1$.  Then the
linear system $|D + t1|$ is of type $g_5^2$ and is different from
$K_C-t$, which gives a contradiction. It follows that $t_i \ne a$. The
image of $\pi$ under the projection $pr_l: \PP^4 \to \PP^2$ is a line
in $\PP^2$ containing four points of the image of $D$. Let $p$ be the
fifth intersection point of this line with $C'$. It follows that $|D|
= |L_p|$, and the statement of the Lemma follows.

   Case 3: $l \cap \pi = \emptyset$. In this case one can repeat the
construction of the case 2, starting from the projection
$pr_l$. E.o.p.

\begin{figure}[htbp]
\begin{center}
   \includegraphics[width=6cm,angle=0]{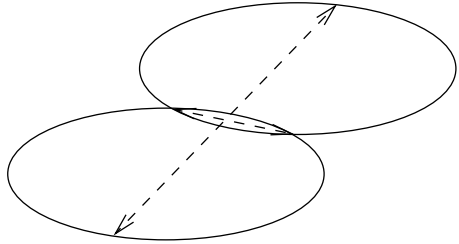}%
% n6
   \caption{Geometry of $W_4^1(C)$.}
   \label{figure:g=5.trigonal.geometry-of-W_4^1}
\end{center}
\end{figure}

It follows that the variety $W = W_4^1(C)$ has two irreducible
components, $W_{(1)}$ and $W_{(2)}$, each isomorphic to the curve $C$
itself: the isomorphism $C \to W_{(1)}$ takes $p$ to $L_p$, and the
isomorphism $C \to W_{(2)}$ takes $p$ to $M_p$. All the points on the
first component, except from the images of $q_1$ and $q_2$,
parametrize globally generated line bundles, while all the points on
the second component parametrize linear systems with a base point. The
Serre duality interchanges the two components (as well as points $q_1$
and $q_2$), as on the Figure
~\ref{figure:g=5.trigonal.geometry-of-W_4^1}.

In the case $C'$ has a cusp $q'$, the arguments of the previous case
can be applied with $q_1 = q_2$. Let $q$ be the preimage of $q'$ on
$C$.  Then the variety $W_4^1(C)$ has two irreducible components, each
isomorphic to $C$, touching each other at the point $q$.

% ------------------------------------------------------------

\subsubsection{Example of a Brill-Noether special K3 of genus 5}
\label{example:trigonal-K3-of-genus-5}
\newcommand{\GG}{{\mathbb G}}
  In this example we construct a polarized K3 surface $(S,a)$ of genus
five such that every hyperplane section of $S$ is trigonal. We will
use the existence of such a K3 in the examples ~\ref{example:g=5,d=5}
and ~\ref{example:g=5,d=4}.

  Let $F=F(0,1,2)$ be a 3-dimensional scroll which can be constructed
as a factor of the affine space $\AA^6$ with coordinates
$(t_0,t_1,x_0,x_1,x_2)$ by the action of the group $\GG_m \times
\GG_m$, where the first factor acts with weights $(1,1,0,-1,-2)$ and
the second one with weights $(0,0,1,1,1)$. The map $\pi: F \to \PP^1$,
where $\phi = t_0/t_1$, gives $F$ a structure of a two-dimensional
projective bundle over $\PP_1$.

  Let $D$ be the subvariety in $F$ defined by the equation $x_0 = 0$.
$D$ has a natural structure of a ruled surface isomorphic to the
two-dimensional scroll $F(1,2)$.  Let $M \in \Pic F$ be the image of
$D$ in the Picard group of $F$, and let $L \in \Pic F$ be the class of
a fiber of $\pi$. Let $\PP = \PP^5$. The map $\phi_M: F \to \PP $
associated with the linear system $|M|$ can be written as $\phi_M =
(x_0: t_0x_1: t_1 x_1: t_0^2x_2: t_0t_1x_2: t_1^2x_2)$.  Note that
$M^3 = 3$ and $M^2L = 1$ in the Chow ring of $F$.

\begin{figure}[htbp]
\begin{center}
\begin{picture}(0,0)%
\includegraphics{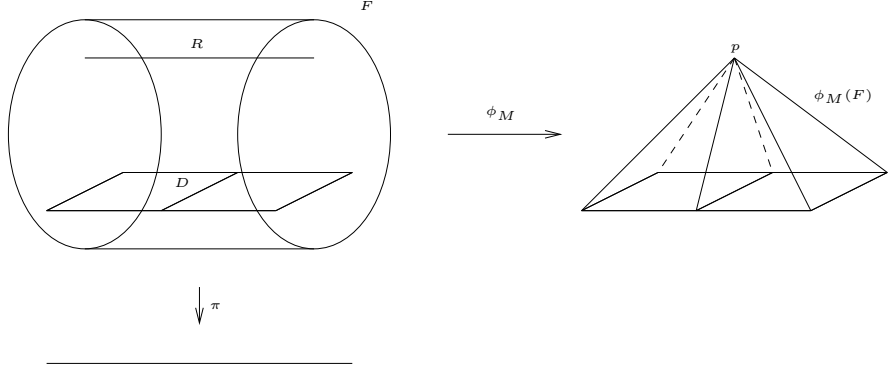}%
% n7
\end{picture}%
\setlength{\unitlength}{3158sp}%
\begingroup\makeatletter\ifx\SetFigFont\undefined%
\gdef\SetFigFont#1#2#3#4#5{%
  \reset@font\fontsize{#1}{#2pt}%
  \fontfamily{#3}\fontseries{#4}\fontshape{#5}%
  \selectfont}%
\fi\endgroup%
\begin{picture}(6920,2895)(593,-2173)
\put(6271,299){\makebox(0,0)[lb]{\smash{\SetFigFont{5}{6.0}{\familydefault}{\mddefault}{\updefault}{\color[rgb]{0,0,0}$p$}%
}}}
\put(4351,-211){\makebox(0,0)[lb]{\smash{\SetFigFont{5}{6.0}{\familydefault}{\mddefault}{\updefault}{\color[rgb]{0,0,0}$\phi_M$}%
}}}
\put(2183,-1733){\makebox(0,0)[lb]{\smash{\SetFigFont{5}{6.0}{\familydefault}{\mddefault}{\updefault}{\color[rgb]{0,0,0}$\pi$}%
}}}
\put(3361,614){\makebox(0,0)[lb]{\smash{\SetFigFont{5}{6.0}{\familydefault}{\mddefault}{\updefault}{\color[rgb]{0,0,0}$F$}%
}}}
\put(2026,314){\makebox(0,0)[lb]{\smash{\SetFigFont{5}{6.0}{\familydefault}{\mddefault}{\updefault}{\color[rgb]{0,0,0}$R$}%
}}}
\put(1906,-774){\makebox(0,0)[lb]{\smash{\SetFigFont{5}{6.0}{\familydefault}{\mddefault}{\updefault}{\color[rgb]{0,0,0}$D$}%
}}}
\put(6923,-84){\makebox(0,0)[lb]{\smash{\SetFigFont{5}{6.0}{\familydefault}{\mddefault}{\updefault}{\color[rgb]{0,0,0}$\phi_M(F)$}%
}}}
\end{picture}

   \caption{Construction of a ``trigonal'' K3, step 1}
   \label{figure:example.g=5.geometry-of-phi_M}
\end{center}
\end{figure}

Let $R$ be the subvariety in $F$ given by the equations $x_1 = x_2 =
0$.  $R$ is a smooth rational curve in $F$ disjoint from $D$ which is
contracted to a point $p = (1,0,0,0,0,0) \in \PP$ under the map
$\phi_M$. One can check that $\phi_M|_{F - R}$ is an embedding. It
follows that $\phi_M(F)$ is a cone over a smooth scrollar surface
$\phi_M(D)$ in $\PP$. The geometry of the map $\phi_M$ is illustrated
by the Figure ~\ref{figure:example.g=5.geometry-of-phi_M}.
% PostScript + LaTeX
% 80% zoom

   Let $H$ be the hyperplane in $\PP$ defined by the vanishing of the
first coordinate; we have $D = \phi_M(F) \intersect H$.

   One can check that $K_F = L - 3M$. Let $S$ be a smooth divisor in
the linear system $|-K_F| = |-L + 3M|$ on $F$. For example, one can
choose $S$ to be given by the ``almost Weierstrass'' equation
$$
   x_0^2x_1 = c_5(t) x_2^3 + c_3(t) x_2 x_1^2 + c_2(t) x_1^3,
$$% 
where $c_k(t) = c_k(t_1,t_2)$ is a generic homogeneous polynomials of
degree $k$. Since $q_S = 0$ and $K_S = 0$, $S$ is a K3 surface. 

   Since $$R \cdot S = (M-L)(M-2L)(3M-L) = -1$$ in the Chow ring of
$F$, we have $R \subset S$. (It is also possible to see it directly,
since the weight consideration implies that all the monomials in the
equation of $S$ must be divisible by $x_1$ or by $x_2$.) Computing the
Chern class gives $N_{S/R} \simeq \OO_R(-2)$, i.e., the map $\phi_M$
contracts the rational curve $R$ on $S$ into a double point $p \in
\phi_M(S)$.

 Let $C = S \intersect D$. We get a diagram
$$
   \xymatrix
   {
      p \xyArrow{\in}{r} & {\PP} \xyArrow{\supset}{r} & H \\
      R \ar@{|->}[u] \xyArrow{\subset}{r} & %
          F \xyArrow{\supset}{r} \ar[u]_{\phi_M} &  %
          D \ar[u]_{\phi_M|_D} \\
      R \xyArrow{\parallel}{u} \xyArrow{\subset}{r} & %
          S \xyArrow{\supset}{r} \xyArrow{\cup}{u} \ar[d]_{\pi|_S} & %
          C \xyArrow{\cup}{u} \ar[ld]^{\pi |_C} \\
      {} & {\PP^1} & {}
   }
$$

Note that
$$ \deg_M F = \deg_M D = M^3 = 3 $$ and
$$ \deg_M S = \deg_M C = M^2 (3M-L) = 8 $$
One can see that $\OO_F(D)|_C \simeq K_C$. It follows that $\phi_M: C
\to H \subset \PP^5$ is the canonical embedding of $C$; in particular,
$C$ has genus five.

The line bundle $\OO_F(M)$, restricted to a fiber $F_t$, is isomorphic
to $\OO_{\PP^2}(1)$, which implies that $\phi_M$ maps every rational
curve $D_t$ to a trisecant line to $C$.  

\begin{figure}[htbp]
\begin{center}
\begin{picture}(0,0)%
\includegraphics{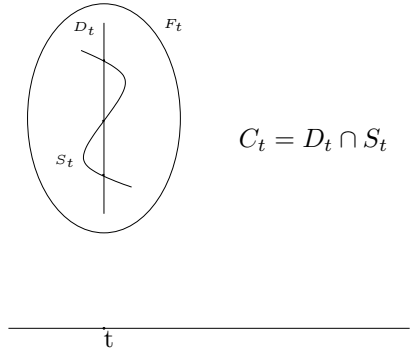}%
% n8
\end{picture}%
\setlength{\unitlength}{3158sp}%
\begingroup\makeatletter\ifx\SetFigFont\undefined%
\gdef\SetFigFont#1#2#3#4#5{%
  \reset@font\fontsize{#1}{#2pt}%
  \fontfamily{#3}\fontseries{#4}\fontshape{#5}%
  \selectfont}%
\fi\endgroup%
\begin{picture}(3174,2706)(439,-2161)
\put(2261,-586){\makebox(0,0)[lb]{\smash{\SetFigFont{10}{12.0}{\familydefault}{\mddefault}{\updefault}{\color[rgb]{0,0,0}$C_t = D_t \cap S_t$}%
}}}
\put(1201,-2161){\makebox(0,0)[lb]{\smash{\SetFigFont{10}{12.0}{\familydefault}{\mddefault}{\updefault}{\color[rgb]{0,0,0}t}%
}}}
\put(816,-716){\makebox(0,0)[lb]{\smash{\SetFigFont{5}{6.0}{\familydefault}{\mddefault}{\updefault}{\color[rgb]{0,0,0}$S_t$}%
}}}
\put(956,324){\makebox(0,0)[lb]{\smash{\SetFigFont{5}{6.0}{\familydefault}{\mddefault}{\updefault}{\color[rgb]{0,0,0}$D_t$}%
}}}
\put(1671,349){\makebox(0,0)[lb]{\smash{\SetFigFont{5}{6.0}{\familydefault}{\mddefault}{\updefault}{\color[rgb]{0,0,0}$F_t$}%
}}}
\end{picture}

   \caption{Construction of a ``trigonal'' K3, step 2}
   \label{picture:genus-5.trigonal}
\end{center}
\end{figure}

Fiber-wise with respect to $\pi$ the geometry of $S$ is shown on the
Figure ~\ref{picture:genus-5.trigonal}.
% PostScript + LaTeX
% 80% zoom

 (Going in the opposite direction, consider the image of a
non-hyperelliptic genus 5 curve $C$ under the canonical embedding
$\phi_{K_C}: C \to \PP^4$. The Riemann-Roch theorem implies that there
is a two-dimensional linear system of quadrics through $\phi_{K_C}(C)$
in $\PP^4$. If $C$ is not trigonal, then $C$ is equal to the
intersection of quadrics trough $C$ by Enriques-Petri theorem; if $C$
is trigonal, intersection of quadrics of through $C$ is a linear
scrollar surface $X$ isomorphic to $F(1,2)$ such that every line on
$X$ is trisecant to $C$.)

I learned about the geometry of scrollar varieties from Miles Reid's
book ``Chapters on Algebraic Surfaces.''

We will use the existence of such a K3 surface in the next two
examples.

%-------------------------------------------------------------------
\subsubsection{Example: $g=5, d=5$}
\label{example:g=5,d=5}
Let $(S,a)$ be a polarized K3 surface of genus $g = 5$, and let $d =
5$. Note that $M(0,a,\beta) = M(0,a,-1) \subset \Pic^{3}(|a|)$.  The
formula ~\ref{eq:virtual-dimension} gives the following geography for
the Brill-Noether loci:
%------------
\begin{comment}
$$
\begin{tabular}{c|cccccc}
    r &    &   &   &   &  &  \\
    2 &    &   &   & 6 & 2 &  \\
    1 &    &10 & 8 & 4 &  &  \\ 
    0 & 10 & 8 & 4 &   &  &  \\
   \hline
      &  0 & 1 & 2 & 3 & 4 & k
\end{tabular}
$$
or, symbolically, 
\end{comment}
% ----------------
$$
\begin{tabular}{c|cccccc}
    r &   &   &   &   &   &  \\
    2 &   &   &   & * & ! &  \\
    1 &   & * & * & ! &   &  \\
    0 & * & * & ! &   &   &  \\
   \hline
      & 0 & 1 & 2 & 3 & 4 & k\\
\end{tabular}
$$
There are three integer points in the domain $D_2$. The variety
$V_3^1(|a|)$ parametrizes trigonal structures on smooth curves in the
linear system $|a|$. The variety $\Hilb^5_{(2,a)}$ parametrizes degree
5 0-subschemes $\xi$ on $S$ which span a two-dimensional subspace in
$\PP^5$. For such a 0-subscheme $\xi$ there is a net of hyperplane
sections of $S$ through $\xi$, each smooth member $C$ of which is a
genus five curve. On such a curve $C$ $\xi$ gives a linear system of
type $g_5^2$. (Note that the Serre dual to $g_5^2$ is $g_3^1$ for
genus 5).

Taking into the account the results of ~\cite{saint-donat}, there are
two cases to consider:

   {\bf Case 1:} There are no trigonal curves in the linear system
$|a|^s$. In this case every smooth curve in $|a|$ is Brill-Noether
general, and it follows that all the Brill-Noether varieties in the
domain $D_2$ are empty. 
\begin{comment}
The geography of the Brill-Noether loci in
this case is
%
$$
\begin{tabular}{c|cccccr}
    r &   &   &   &   &  & \\
    2 &   &   &   & * &  & \\
    1 &   & * & * &   &  & \\
    0 & * & * &   &   &  & \\
   \hline
      & 0 & 1 & 2 & 3 & 4 & k 
\end{tabular}
$$
\end{comment}
%
This is the case of a generic polarized K3 of genus 5. The results of
Saint-Donat imply that such a K3 surface $S$ is an intersection of 3
quadrics in $\PP^5$.

{\bf Case 2:} All smooth curves in the linear system $|a|$ are
trigonal.  Since every smooth trigonal genus 5 curve has precisely one
trigonal structure, the varieties $V_3^1(|a|)$ and
$\Hilb^5_{(2,a)}(S)$ have dimension five, which is greater than one
predicted by the formula ~\ref{eq:virtual-dimension}). 
\begin{comment}
It follows that the geography of the Brill-Noether
varieties in this case is 
%      - this is not true; we do not know anything about rank 2
%
$$
\begin{tabular}{c|ccccc}
    2 &   &   &   & * &    \\
    1 &   & * & * & * &    \\
    0 & * & * & * &   &    \\
   \hline
      & 0 & 1 & 2 & 3 & k 
\end{tabular}
$$
%
\end{comment}

Remark 1: results of Saint-Donat (~\cite{saint-donat}) imply that such
a K3 surface $S$ we have $\rank \Pic S \ge 2$, and $\Pic S$ should
contain the sublattice
$$
  \begin{tabular}{|cc|}
      8 & 3 \\
      3 & 0 
  \end{tabular}
$$ (which is obviously true for the K3 constructed in the previous
example).
%
\begin{comment}
%    i.e., $\Pic S$ contains $Z C + Z E$, where $E$ is an elliptic curve; 
%    the linear system $|E|$ restricts to a $g_3^1$ on $C$. 
%    If $H = a C + b E$ ($H$ is the stability parameter), then
%    $(H,C) \ZZ = (8a + 3 b) \ZZ = \ZZ$, since g.c.d. = 1
%            $(H,C) = 1: 8a + 3b = 1$
%            is possible only if $a < 0$ , $b > 0$, or
%                               $a > 0$, $b < 0$, - 
%                               but such an $H$ is not ample.
\end{comment}
Note that $\Pic S$ does not contain a nef divisor $H$ such that $(a,H)
= 1$. It follows that the main theorem of this text does not imply
that the variety $BN_4(2,a,-1) \in D_2$ is empty, since both of the
the stability criteria we proved are not applicable in this case.

Remark 2: Let $(C,B) \in V_3^1(|a|)$. Note that both $B$ and $K_C
B^{-1}$ are globally generated. By Lemma
~\ref{generic-extensions-are-acceptable}, there is a unique up to an
isomorphism locally free globally generated extension $E$ of $(C,B)$
of rank 2. Since $\rho(C,B) < 0$ (i.e., the generic genus 5 curve is
not trigonal), the Riemann-Roch formula implies that $E$ is not simple
and, in particular, is not stable.
% (compare Example
%~\ref{example:case-of-a-negative-rho}).

Remark 3: Let $\xi \in \Hilb^5_{(2,a)}(S)$ be a simple 0-cycle.  We
claim that $\xi$ is Caley-Bacharash with respect to $|a|$.  First
proof: $\xi$ induces a linear system of type $g_5^2$ on every curve $C
\in |a|$ through $\xi$. Let $C'$ be a plane model of $C$ given by this
linear system. Then $\xi$ is the intersection of $C'$ with a line $l$
on $\PP^2$. Note that $K_C = \OO_{\PP^2}(2h)|_C$.  Since every quadric
on $\PP^2$ through $\xi$ is of the form $l + l'$, where $l' \in |h|$,
$\xi$ is Caley-Bacharash with respect to $|a|$. Second proof: if $\xi$
is not Caley-Bacharash, then it contains a sub-cycle $\xi'$ of degree
4 and speciality index 2, which contradicts to the Martens-Mumford
theorem applied to any curve through $\xi'$.  

It follows that the extension defined in the Lemma%
~\ref{generic-extensions-are-acceptable} %% 1 --> r
is locally free. (We have already proved that it is not stable.)

%---------------------------------------------------------------
%
\subsubsection{Example: $g=5, d=4$}
\label{example:g=5,d=4}
Let $(S,a)$ be a polarized K3 surface of genus $g = 5$, and let $d =
4$. Note that $M(0,a,\beta) = M(0,a,0) \subset \Pic^{4}(|a|)$.  The
formula ~\ref{eq:virtual-dimension} gives the following table of the
expected dimensions of the Brill-Noether loci: 
$$
\begin{tabular}{c|cccccc}
      r &    &   &   &   &    &   \\
      2 &    &   &   &   &  2 &   \\
      1 &    &   & 8 & 5 &  0 &   \\
      0 &    & 9 & 6 & 1 &    &   \\
   \hline
        & 0  & 1 & 2 & 3 &  4 & r    
  \end{tabular}
$$
or, symbolically,
$$
\begin{tabular}{c|cccccc}
      r &   &   &   &   &   &     \\
      2 &   &   &   &   & * &     \\
      1 &   &   & * & * & ! &     \\
      0 &   & * & * & ! &   &     \\
   \hline
        & 0 & 1 & 2 & 3 & 4 & k  
\end{tabular}
$$
The variety $V_4^2(|a|^s)$ is empty by the Martens-Mumford theorem,
and the geometry of $V_4^1(|a|^s)$ depends on the geometry of curves
in the linear system $|a|$.
 
    {\bf Case 1:} There are no trigonal curves in $|a|^s$. If $C$ is a
Brill-Noether general curve of genus 5, it is easy to prove that $V(C)
:= V_4^1(C)$ is a smooth irreducible curve of genus 9. All linear
systems parametrized by $V(C)$ are globally generated, and the Serre
duality acts as an involution on $V(C)$.  There are many more results
on the geometry of $V(C)$ in ~\cite{tyurin-o-peresechenii-kvadrik}.

  Let $V(|a|) = V_4^1(|a|^s)$. This variety is fibered over $|a|^s$
with fibers $V_4^1(C)$.

The moduli space $S':= M_H(2,a,0)$ is two-dimensional, and thus is a
K3-surface. There is a well-defined morphism $E: V(|a|) \to S'$.
$$
\begin{tabular}{c|ccccc}
      2 &    &          &   & $S'$ & \\
      1 &    & *        & * &      & \\
      0 &  * & $V(|a|)$ &   &      & \\
   \hline
        &  1 & 2        & 3 & 4    & k   
\end{tabular}
$$
{\bf Case 2:} All curves in the linear system $|a|^s$ are
trigonal. The example of such a surface was constructed in the Example
~\ref{example:trigonal-K3-of-genus-5}.

   Let $C$ be a trigonal curve of genus 5. We proved in the Example
~\ref{example:two-components-of-W_4^1(C)} that $V(C):=W_4^1(C)$ is
one-dimensional and has two irreducible components, $V_{(1)}(C)$ and
$V_{(2)}(C)$, intersecting each other at two points or one double
point. Generic points on the first component parametrize globally
generated line bundles on $C$, and all points on the second component
parametrize linear systems with base points. The Serre duality
interchanges these two components.  It follows that $V(|a|)=
V_4^1(|a|^s)$ has two irreducible components, $V_{(1)}(|a|)$ and
$V_{(2)}(|a|)$. Note that the generic point of $V_{(1)}(|a|)$ is not
in the image of the correspondence $\cA$.

%  v odnoj komponente lezhat B - ne gg, zato A - gg i mozho rashirit' do 
%  vektornogo rassloeniya. Ono tochno ne budet globally generated. 
%  rho = 1 => ono mozhet byt' stabil'nym.
%  Zametim, chto C - tochno ne obrazuyuschaya Pikara, i (C,H) = 1
%  tozhe nevozmozhno podobrat'.  hotelos' by razobrat' rukami!!!!

{\bf Lemma:} $S'= M(2,a,0)$ does not contain any globally generated
vector bundles.

Proof: Assume that $S'$ contains a globally generated vector bundle
$E$. Note that the Grassmanian variety $\Gr(2,h^0(S,E))$ has positive
dimension.  Let us choose a generic point $V \in \Gr(2,h^0(S,E))$.  By
lemma %
~\ref{generic-factors-in-globally-generated-case-are-acceptable}, %
the evaluation map $e_{(V,E)}$ has a good degeneration, and its
cokernel is a point $v = (C,B) \in V_4^1(|a|^s)$ such that both $B$
and $K_C B^{-1}$ are globally generated. Since $V(|a|)$ does not
contain such a point $p$, we get a contradiction.

{\bf Corollary} Variety $\cA_{2,0}^{4,2}$ is empty.

{\bf Remark.}  \newcommand{\Coh}{{\rm Coh}}
One can still define the extension morphism
$$ E: V(|a|) \to \Coh_S$$, where $\Coh_S$ is the stack of coherent
sheaves on S. The stack $\Coh_S$ has a substack $\Coh_{(1)}$ of
non-locally free sheaves and a substack $\Coh_{(2)}$ of non-globally
generated sheaves. Under $E$, $V_{(i)}(|a|)$ is the preimage of
$\Coh_{(i)}$.

%-------------------------------------------------------------------------
%

%
%----------------------------------------------------------------------

\subsubsection{Example of a Brill-Noether special K3 of genus 6}
\label{example:K3-of-genus-6-with-g_5^2}

   In this example we construct a polarized K3-surface $(S,a)$ of
genus six such that every hyperplane section of $S$ is isomorphic to a
plane quintic. We will use the results of this Example in the Example
~\ref{g=6,d=6}.

   Let $(S,B)$ be a polarized genus K3 surface of genus two. Thus $B
\in \Pic(S)$, $B^2 = 2$, and the map $\phi_B: S \to |B|\dual \simeq
\PP^2$ realizes $S$ as a two-sheeted cover of $\PP^2$ branched along a
curve $R \subset \PP^2$ of degree six. For a degree six curve $R$
given by the equation $f(x_0,x_1,x_2) = 0$ the surface $S$ can be
realized as a subscheme in the total space of the line bundle
$\OO_{\PP^2}(3)$ given by the equation
$$
    \zeta^2 = f(x_0,x_1,x_2)
$$

   Let $R$ be a plane sextic which has a tritangent line $l$, i.e.,
$f$, restricted to $l$, is a square of a homogeneous form of degree
$3$. For example, we can take
$$ f = x_1^6 + x_2^6 + x_0^5 x_2, l = (x_2 = 0),
$$ This is a very degenerate case, since all three points at which $l$
is bitangent to $R$ coincide. We have $\phi_B^{-1}(l) = \Gamma +
\Gamma'$, where $\Gamma$ and $\Gamma'$ are smooth rational curves, and
$(B,\Gamma) = 1$.

Note that the linear system $a=2B+\Gamma$ has genus 6, and
$(2B+\Gamma,B)_S = 5$. It follows that $|B|$, restricted to a curve in
$|2B + \Gamma|$, gives a linear system of type $g_5^2$.  In
particular, the image of $\phi_{2B + \Gamma}: S \to |2B+\Gamma|\dual
\simeq \PP^6$ is as a projective K3 surface of genus such that all of
its hyperplane sections are isomorphic to plane quintics (For the
non-constructive approach, see ~\cite{saint-donat}).

  Note that since $(2B + \Gamma, \Gamma) = 0$, the map $\phi_{2B +
\Gamma}$ contracts $\Gamma$ into a point $p_{\Gamma}$.

Let $V_1 = H^0(S,B)$ and $V_2 = H^0(S,2B)$. Note that the
canonical map
$$
   \mathop{\rm Sym}\nolimits^2 V_1 \to V_2
$$ is an isomorphism. Considering the adjunction sequences
$$ 0 \to \OO(2B) \to \OO(2B + \Gamma) \to \OO(2B + \Gamma)|_{\Gamma}
   \to 0
$$
and
$$ 0 \to \OO(B) \to \OO(2B) \to \OO(2B)|_{B} \to 0
$$ on $S$, one can construct the diagram
%  to consider the square
%
%                2B + Gamma 
%                   |
%            B  ->  2B
%
$$
  \xymatrix
  {
    S \ar[r]^{ \phi_{2B + \Gamma} } 
       \ar[rd]^{\phi_{2B}}           
       \ar[d]_{\phi_{B}}         &       
    {\PP^6} \ar@{-->}[d]^{pr_{\Gamma}} \\
    {\PP^2} \ar@<1mm>[r]^{2h}         & 
    {\PP^5} \ar@{-->}@<1mm>[l]^{pr_{B_0}}
  }
$$ where $pr_{\Gamma}$ is a projection map with center $p_{\Gamma}$
and $pr_{B_0}$ is a projection from the (2-dimensional) linear span of
any curve $B_0 \in |B|$.

%   This gives an example of a polarized K3 surface $(S,2B+\Gamma)$ of
%genus 6 for which every hyperplane section is isomorphic to a plane
%quintic, which we will need in the following example.

%----------------------------------------------------------------------

\subsubsection{Example: $g=6,d=6$}
\label{g=6,d=6}
Let $(S,a)$ be a polarized K3 surface of genus $g = 6$, and let $d =
6$. Note that $M(0,a,\beta) = M(0,a,-1) \subset \Pic^{4}(|a|)$.  The
formula ~\ref{eq:virtual-dimension} gives the following expected
dimensions of the Brill-Noether loci:
$$
\begin{tabular}{c|ccccccc}
      r &    &    &   &   &   &    & \\
      3 &    &    &   &   &   & 0  & \\
      2 &    &    &    & 8 & 4 &   & \\
      1 &    & 12 & 10 & 6 & 0 &   & \\
      0 & 12 & 10 &  6 & 0 &   &   & \\
   \hline
        & 0  & 1  & 2  & 3 & 4 & 5 & k
\end{tabular}
$$
%
%
\begin{comment}
or, symbolically, 
%
$$
\begin{tabular}{c|ccccccc}
      r &   &   &   &   &   &   & \\
      3 &   &   &   &   &   & * & \\
      2 &   &   &   & * & * &   & \\
      1 &   & * & * & * & ! &   & \\
      0 & * & * & * & ! &   &   & \\
   \hline
        & 0 & 1 & 2 & 3 & 4 & 5 & k 
\end{tabular}
$$
%
%
\end{comment}
% 
The Martens-Mumford theorem implies that $V_4^2(C)$ is empty for any
genus 6 curve. It follows that $\Hilb^6_{3,a}(S)$ is empty. Let
$V(|a|)=V_4^1(|a|)$ and $M= M_H(3,a,-1)$. Symbolically, the geography
of the Brill-Noether loci is
$$
\begin{tabular}{c|ccccccc}
   3 &   &   &   &   &   & $M$  \\
   2 &   &   &   & * & * &      \\
   1 &   & * & * & * &   &      \\
   0 & * & * & $V(|a|)$ &&&     \\
   \hline
     & 0 & 1 & 2 & 3 & 4 & 5   
\end{tabular}
$$
For a genus 6 curve one has $\rho(g_4^1) = 0$ and $\dim V_4^1(C) \le
1$ by the Martens-Mumford theorem. We will consider the following two
cases:
 {\bf Case 1:} There is a smooth Brill-Noether general curve in the
linear system $|a|$. Let $V(C)= V_4^1(C)$.  It is well-known (cf.
~\cite{ACGH}) that for a Brill-Noether general curve $C$ of genus 6
the variety $V(C)$ is a union of 5 distinct points.
%
%, and thus $V(|a|)$ is finite over $|a|^s$.
% (is it true?)
%
  One can easily check that every point in $V_4^1(C)$ parametrizes
such a line bundle $L$ on $C$ that both $L$ and $K_C L^{-1}$ are
globally generated. 

The moduli space $M$ has dimension 0. The results of ~\cite{mukai-84}
imply that $M$ is either empty or a reduced irreducible variety of
dimension 0 (i.e., is a point).  If $\Pic S = \ZZ a$, there is a
well-defined map $E: V(|a|) \to M$, which proves that $M$ is not
empty, and since $E$ realizes $V(|a|)$ as a Grassmanian fibration over
$M$, it follows that $V(|a|)$ is birational to the Grassmanian variety
$\Gr(3,5)$ of dimension 6. 
%
%In particular, $V(|a|)$ is irreducible. 
%
Note that there exist a degree 5 morphism $V(|a|) \to |a|^s \subset
|a| \simeq \PP^6$.

%     kak monodromiya perestavlyaet linejnye sistemy tipa 1 i tipa 2

%  (ACGH p. 210),

%  Vopros: V - svyazno, ili net?  (t.e. monodromoya - tranzitivna?)
%  pohozhe, chto V imeet 2 svyaznye komponenty
%          (g_4^1, pysekaemye pryamymi, i vysekaemye konikami)

%     r= 3 &           M 
%     r= 2 &       * * 
%     r= 1 &   * * *  
%     r= 0 & * * V  
%     \hline
%     k =   0 1 2 3 4 5 

{\bf Case 2:} All curves in the linear system $|a|$ have $g_5^2$.  An
example of such a surface was constructed in
~\ref{example:K3-of-genus-6-with-g_5^2}.

  Let $C$ be a genus six curve which can be realized as a (smooth)
plane quintic, and let $|h|$ be the corresponding linear system of
type $g_5^2$.  Projection map with center $p \in C$ gives a linear
system of type $g_4^1$ on $C$, and one can prove that every $g_4^1$ on
$C$ arises in such a way. It follows that $V_4^1(C) \simeq C$. Every
linear system $B \in V_4^1(C)$ is globally generated, and every $K_C -
B = 2h - (h-p) = h+p$ has a base point $p$. It follows that no $B \in
V(C)$ have locally free extensions.

   Vice versa, the only stable sheaf $E \in M$ is locally free by the
(~\cite{mukai-84}, proposition 3.3), but one can see as in the
previous examples that $E$ is not globally generated. It follows that
the variety $\cA_{3,0}^{5,3}$ is empty.

% -- rashirenie lyubogo B ne loka'no svobodno --
%
%  F:  V --> M
%      6     0
%
%      F ne opredelen: oni mogut byt' stabil'ny
%      (hotya ne yasno, kak eto garantirovat': Pic > Z), no ne locally free.
%
%      naoborot, rassmotrim E. Esli E gg, to u nego est' (C,B),
%      B = g_4^1 - gg, A = g_6^2 - gg - tak ne byvaet =>
%      E - *ne* global'no porozhdeno.
%

%

\end{document}